\documentclass[10pt,a4paper,equation]{article}

\usepackage[margin=7em]{geometry}

\usepackage[greek,american]{babel}
\usepackage[iso-8859-7]{inputenc}
\usepackage{amssymb, amsfonts}
\usepackage[dvips]{graphicx}
\usepackage{amsmath, amsthm}
\usepackage{epsfig}
\usepackage{latexsym}
\usepackage{makeidx}

\newtheorem{theorem}{Theorem}[section]

\numberwithin{equation}{section}

\newcommand{\EE}{\textbf{E}}
\newcommand{\PP}{\mbox{\bf P}}
\newcommand{\N}{{\mathbb N}}
\newcommand{\R}{{\mathbb R}}
\newcommand{\Np}{\mathbb{N}^+}

\newcommand{\gl}{\lambda}

\newcommand{\gep}{\varepsilon}

\newcommand{\sm}{{\raise0.3ex\hbox{$\scriptstyle \setminus$}}}

\newcommand{\beq}{\begin{equation}}
\newcommand{\eeq}{\end{equation}}

\makeatletter
\def\blfootnote{\gdef\@thefnmark{}\@footnotetext}
\makeatother

\title{Functional limit theorems for the P\'olya and $q$-P\'olya urns}

\author{
	Dimitris Cheliotis\thanks{Department of Mathematics,
	Panepistimiopolis, Athens 15784, Greece.
        Email: {\tt dcheliotis@math.uoa.gr}.
	}
\and Dimitra Kouloumpou \thanks{Techological Educational Institute of Pelloponnese, Antikalamos 24100, Kalamata, Greece.  Email: {\tt dimkouloubou@teipel.gr}.}
	}

\begin{document}
%\date{}
\maketitle



\begin{abstract} For the plain P\'olya urn with two colors, black and white, we prove a functional central limit theorem for the number of white balls assuming that the initial number of black balls is large. Depending on the initial number of white balls, the limit is either a pure birth process or a diffusion. We also prove analogous results for the $q$-P\'olya urn, which is an urn where, when picking a ball, the balls of one color have priority over those of the other.  
\end{abstract}

\section{Introduction and results} 

\subsection{The models} \label{Models}

\hspace{2.5ex} \textbf{The P\'olya urn.} This is the model where in an urn that has initially $r$ white and $s$ black balls we draw, successively, uniformly, and at random, a ball from it and then we return the ball back together with $k$ balls of the same color as the one drawn. The number $k\in\Np$ is fixed. Call $A_n$ and  $B_n$ the number of white and black balls respectively after $n$ drawings. The most notable result regarding its asymptotic behavior is that the proportion of white balls in the urn after $n$ drawings, $A_n/(A_n+B_n)$, converges almost surely as $n\to\infty$ to a random variable with distribution Beta$(r/k, s/k)$. Our aim in this work is to examine whether the entire path $(A_n)_{n\ge 0}$ after appropriate natural transformations converges to a stochastic process.

Standard references for the theory and the applications of P\'olya urn and related models are \cite{JoKo} and \cite{Mah}.

\textbf{The $q$-P\'olya urn}. This is a $q$-analog of the P\'olya urn (see \cite{Gas}, \cite{KacCheung} for more on $q$-analogs) introduced in \cite{Kup} and studied further in \cite{Char12} (see also \cite{Char16}). A $q$-analog of a mathematical object $A$ is another object $A(q)$ so that when $q\to1$, $A(q)$ ``tends'' to $A$. Take $q\in(0, \infty)\sm\{1\}$. The $q$-analog of any $x\in \mathbb{C}$ is defined as
\begin{equation} \label{qNumber}
[x]_q:=\frac{q^x-1}{q-1}.
\end{equation}
Note that $\lim_{q\to1} [x]_q=x$. Now consider an urn that has initially $r$ white and $s$ black balls, where $r, s\in\N, r+s>0$. We perform a sequence of additions of balls in the urn according to the following rule. If at a given time the urn contains $w$ white and $b$ black balls ($w, b\in \N, w+b>0$),  then we add $k$ white balls with probability
\begin{align} \label{qPWhite}  \PP_q(\text{white})&=\frac{[w]_q}{[w+b]_q}. \\
\intertext{Otherwise, we add $k$ black balls, and this has probability}
\PP_q(\text{black})&=1-\PP_q(\text{white})=q^w \frac{[b]_q}{[w+b]_q}. \label{qPBlack}
\end{align}    

To understand how the $q$-P\'olya urn works, it helps to realize the probabilities $\PP_q(\text{white}), \PP_q(\text{black})$ through a natural experiment. 

If $q\in (0, 1)$, then we put the balls in a line with the $w$ white coming first and the $b$ black following.   To pick a ball, we go through the line, starting from the beginning and picking each ball with probability $1-q$ independently of what happened with the previous balls. If we finish the line without picking a ball, we start from the beginning. Once we pick a ball, we return it to its position together with $k$ balls of the same color. Given these rules, the probability of picking a white ball is
\begin{equation} \label{LinePick}
(1-q^w)\sum_{j=0}^\infty (q^{w+b})^j=\frac{1-q^w}{1-q^{w+b}}=\frac{[w]_q}{[w+b]_q},
\end{equation}
which is \eqref{qPWhite}, because before picking a white ball, we will go through the entire list a random number of times, say $j$, without picking any ball and then, going through the white balls, we pick one (probability $1-q^w$).  

 If $q>1$, we place in the line first the black balls and we go through the list picking each ball with probability $1-q^{-1}$. According to the above computation, the probability of picking a black ball is 
 $$\frac{[b]_{q^{-1}}}{[w+b]_{q^{-1}}}=q^w \frac{[b]_q}{[w+b]_q},$$
which is \eqref{qPBlack}.

We extend the notion of drawing a ball from a $q$-P\'olya urn to the case where exactly one of $w, b$ is infinity. Then the probability to pick a white (resp. black) ball is determined again by \eqref{qPWhite} (resp. \eqref{qPBlack}), where this is understood as the limit of the right hand side as $w$ or $b$ goes to $\infty$. For example, assuming that $w=\infty$ and $b\in \N$, we have $\PP_q(\text{white})=1$ if $q<1$ and $\PP_q(\text{white})=q^{-b}$ if $q>1$. Again these probabilities are realized through the experiment described above. Thus, we can run the process even if we start with an infinite number of balls from one color and finite from the other.

\subsection{P\'olya urn. Scaling limits}  \label{Polya_Scaling_Sec}

For the results of this section, we consider an urn whose initial composition depends on $m\in\Np$. It is $A_0^{(m)}$ and $B_0^{(m)}$ white and black balls respectively. After $n$ drawings, the composition is $A_n^{(m)}, B_n^{(m)}$.

To see a new process arising out of the path of $(A_n^{(m)})_{n\ge0}$ we start with an initial number of balls that tends to infinity as $m\to\infty$. We assume then that $B_0^{(m)}$ grows linearly with $m$. Regarding $A_0^{(m)}$, we study three regimes:

\begin{itemize}
\item[a)] $A_0^{(m)}$ stays fixed with $m$. 
\item[b)] $A_0^{(m)}$ grows to infinity but sublinearly with $m$. 
\item[c)] $A_0^{(m)}$ grows linearly with $m$.  
\end{itemize}

The regime where $A_0^{(m)}$ grows superlinearly with $m$  follows by regime b) by changing the roles of the two colors.

In the regimes a) and b), the scarcity of white balls has as a result that the time between two consecutive drawings of a white ball is of order $m/A_0^{(m)}$ (the probability of picking  a white ball in the first few drawings is approximately $A_0^{(m)}/m$, which is small). We expect then that speeding up time by this factor we will see a birth process. And indeed this is the case as our first two theorems show. 

All processes appearing in this work  with index set $[0, \infty)$ and values in some Euclidean space $\R^d$ are elements of $D_{\R^d}[0, \infty)$, the space of functions $f:[0, \infty)\to\R^d$ that are right continuous and have limits from the left of each point of $[0, \infty)$. This space is endowed with the Skorokhod topology, and convergence in distribution of processes with values on that space is defined through that topology.

We remind the reader that the negative binomial distribution with parameters $\nu\in(0, \infty)$ and $p\in(0, 1)$ is the distribution with support in $\N$ and probability mass function
\begin{equation}\label{NBDensity}   f(x)={x+\nu-1 \choose x}p^\nu(1-p)^x
\end{equation}
for all $x\in\N$. When $\nu\in\Np$, this is the distribution of the number of failures until we see the $\nu$-th success in a sequence of independent trials, each having probability of success $p$. For a random variable $X$ with this distribution, we write $X\sim NB(\nu, p)$.

\begin{theorem} \label{PolyaPathRegime1} Fix $w_0\in\Np$ and $b_0\ge0$. If  $A_0^{(m)}=w_0$ and $\lim_{m\to\infty}B_0^{(m)}/m=b_0$, then the process $(k^{-1}\{A^{(m)}_{[mt]}-A_0^{(m)}\})_{t\ge0}$ converges in distribution, as $m\to\infty$, to an inhomogeneous in time pure birth process $Z=(Z_t)_{t\ge0}$ such that for all $0\le t_1<t_2, j\in\N$, the random variable $Z(t_2)-Z(t_1)|Z(t_1)=j$ has distribution $NB\big(\frac{w_0}{k}+j, \frac{t_1+(b_0/k)}{t_2+(b_0/k)}\big)$. Equivalently, $Z$ has rates $\gl_{t, j}=(k j+w_0)/(k t+b_0)$ for all $(t, j)\in [0, \infty)\times \N$.
\end{theorem}

\begin{theorem} \label{PolyaPathRegime2}  If  $A_0^{(m)}=:g_m$ with $g_m\to\infty, g_m=o(m)$ and $\lim_{m\to\infty}B_0^{(m)}/m=b_0$ with $b_0>0$ constant, then the process $(k^{-1}\{A_{[tm/g_m]}^{(m)}-A_0^{(m)}\})_{t\ge0}$, as $m\to\infty$, converges in distribution to the Poisson process on $[0, \infty)$ with rate $1/b_0$. 
\end{theorem}

Next, we look at regime c), i.e., in the case that at time 0 both black and white balls are of order $m$. In this case, the normalized process of the number of white balls has a non-random limit, which we determine, and then we study the fluctuations of the process around this limit. 

\begin{theorem}   \label{PolyaPathLinear}
Assume that $A_{0}^{(m)}, B_{0}^{(m)}$ are such that $\lim_{m\to\infty} \frac{A_{0}^{(m)}}{m}=a, \frac{B_{0}^{(m)}}{m}=b$ where $a, b\in[0, \infty)$ are not both zero. Then the process $(A_{[mt]}^{(m)}/m)_{t\ge0}$, as $m\to\infty$, converges in distribution to the deterministic process $X_t=\frac{a}{a+b}(a+b+k t), t\ge0$.
\end{theorem}

The limit $X$ is the same as in an urn in which we add at each step $k$ white or black balls with corresponding probabilities
$a/(a+b), b/(a+b)$, that is, irrespective of the composition of the urn at that time. 

To determine the fluctuations of the process  $(A_{[mt]}^{(m)}/m)_{t\ge0}$ around its $m\to\infty$ limit, $X$, we let
$$C_t^{(m)}=\sqrt{m}\bigg( \frac{A_{[mt]}^{(m)}}{m}-X_t\bigg)$$
for all $m\in\Np$ and $t\ge0$.

\begin{theorem} \label{PolyaDiffusionLimit}  Let $a, b\in[0, \infty)$, not both zero, $\theta_1, \theta_2\in \R$, and assume that $A_0^{(m)}:=[am+\theta_1\sqrt{m}], B_0^{(m)}=[bm+\theta_2\sqrt{m}]$ for all large $m\in\N$. Then the process $(C_t^{(m)})_{t\ge0}$ converges in distribution, as $m\to\infty$, to the unique strong solution of
the stochastic differential equation
\begin{align}
Y_0&=\theta_1, \\
dY_t&=\frac{k}{a+b+k t} \bigg\{Y_t-\frac{a}{a+b}(\theta_1+\theta_2)\bigg\} \,dt+k\frac{\sqrt{ab}}{a+b}\,dW_t,
\end{align}
which is
\begin{equation} \label{PolyaDiffLimit}
Y_t=\theta_1+\frac{b\theta_1-a\theta_2}{(a+b)^2 } kt+k\frac{\sqrt{ab}}{a+b} (a+b+k t) \int_0^t \frac{1}{a+b+ks}\, dW_s.
\end{equation}
$W$ is a standard Brownian motion
\end{theorem}

\noindent \textbf{Remark}. Functional central limit theorems for P\'olya type urns have been proven with increasing generality in the works \cite{Go93}, \cite{BaHu02}, \cite{Ja}. The major difference with our results is that in theirs, the initial number of balls, $A_0^{(m)}, B_0^{(m)}$, is fixed. More specifically:

1) Gouet (\cite{Go93}) studies urns with two colors (black and white) in the setting of Bagchi and Pal (\cite{BaPa85}). According to that, when a white ball is drawn, we return it in the urn together with $a$ white and $b$ black balls, while if a black ball is drawn, we return it together with $c$ white and $d$ black. The numbers $a, b, c, d$ are fixed integers (possibly negative), the number of balls added to the urn is fixed (that is $a+b=c+d$), and balls are drawn uniformly form the urn. The plain P\'olya urn is not studied in that work because, according to the author, it has been studied by Heyde in \cite{He77}. However, for the P\'olya urn, \cite{He77} discusses the central limit theorem and the law of the iterated logarithm. In any case, following the techniques of Heyde and Gouet one can prove the following. Assume for simplicity that $k=1$ and let $L=:\lim_{n\to\infty}\frac{A_n}{n}$. The limit exists with probability one because of the martingale convergence theorem. Then 
$$\left\{\sqrt{n}\left(t\frac{A_{n/t}}{n}-L\right)\right\}_{t\ge0} \overset{d}{\to} \{W_{L' (1-L') t}\}_{t\ge0}$$ 
as $n\to\infty$. $W$ is a standard Brownian motion and $L'$ is a random variable independent of $W$ and having the same distribution as $L$. On the other hand, de-Finetti's theorem gives easily the more or less equivalent statement that, as $n\to\infty$,
$$\left\{\sqrt{n}\left(\frac{A_{nt}}{nt}-L\right)\right\}_{t\ge0} \overset{d}{\to}  \{W_{L' (1-L')/t}\}_{t\ge0}$$ 
with $W, L'$ as before.

2) Bai, Hu, and Zhang (\cite{BaHu02}) work again in the setting of Bagchi and Pal, but now the numbers $a, b, c, d$ depend on the order of the drawing and are random. The requirement that each time  we add the same number of balls is relaxed.

3)  Janson (\cite{Ja}) considers urns with many colors, labeled $1, 2, \ldots, l$, where after each  drawing, if we pick a ball of color $i$, we place in the urn balls of every color according to a random vector $(\xi_{i, 1}, \ldots, \xi_{i, l})$ whose distribution depends on $i$ ($\xi_{i, j}$ is the number of balls of color $j$ that we add in the urn). Also, each ball is assigned a certain nonrandom activity that depends only on its color, and then the probability to pick a certain color at a drawing equals the ratio of the total of the activities of all balls of that color to the total of the activities of all balls present in the urn at that time. A restriction in that work is that there is a color $i_0$ so that starting the urn with just one ball and this ball has this color, there is positive probability to see in the future every other color. This excludes the classical P\'olya urn that we study.

\subsection{$q$-P\'olya urn. Basic results} \label{qPolyaBasic}

We recall some notation from $q$-calculus (see \cite{Char16}, \cite{KacCheung}). For $q\in(0, \infty)\sm\{1\}, x\in \mathbb{C}, k\in \Np$, we define
\begin{align}
[x]_q&:=\frac{q^x-1}{q-1} & \text{the $q$-number of } x,\\
[k]_q!&:=[k]_q [k-1]_q \cdots [1]_q   & \text{the $q$-factorial},\\
[x]_{k, q}&:=[x]_q [x-1]_q \cdots [x-k+1]_q  & \text{  the $q$-factorial of order } k,\\
{x\brack k}_{q}&:=\frac{[x]_{k, q}}{[k]_q!}  & \text{ the $q$-binomial coefficient}\\
(x; q)_\infty&:=\prod_{i=0}^\infty(1-xq^i)  \text{ when }q\in[0, 1) & \text{the $q$-Pochhammer symbol} \label{qPoch},
\end{align}
We extend these definitions in the case $k=0$ by letting $[0]_q!=1, [x]_{0, q}=1$. 

Now consider a $q$-P\'olya urn that has initially $r$ white and $s$ black balls, where $r\in\N\cup\{\infty\}$ and $s\in\N$. Call $X_n$ the number of drawings that give white ball in the first $n$ drawings. Its distribution is specified by the following.

\noindent \textbf{Fact 1:}  Let  $a:=r/k$ and $b:=s/k$. 

\noindent (i) If $r\in \N$, then the probability mass function of $X_n$ is
\begin{align} \label{QPolyaWhiteDistr}
\PP\left(X_{n}=x\right)&=q^{k(n-x)(a+x)}  \frac{{-a \brack x}_{q^{-k}} {-b \brack n-x}_{q^{-k}}}{{-a-b\brack n}_{q^{-k}}}=q^{-sx}  \frac{{a+x-1 \brack x}_{q^{-k}} {b+n-x-1\brack n-x}_{q^{-k}}}{{a+b+n-1\brack n}_{q^{-k}}} \\
&=q^{-kx(b+n-x)}  \frac{{-a \brack x}_{q^{k}} {-b \brack n-x}_{q^{k}}}{{-a-b\brack n}_{q^{k}}} \label{QPWDQ}
\end{align}
for all $x\in\N$.

\noindent (ii) If $r=\infty$ and $q>1$, then the probability mass function of $X_n$ is
\begin{align} \label{QPolyaWhiteDistrRInfty}
\PP\left(X_{n}=x\right)&=q^{-sx}(1-q^{-k})^{n-x}  {b+n-x-1\brack n-x}_{q^{-k}} \frac{[n]_{q^{-k}}!}{[x]_{q^{-k}}!} 
\end{align}
for all $x\in\N$.

Relation \eqref{QPolyaWhiteDistr} is (3.1) in \cite{Char12} where it is proved through recursion. In Section \ref{qPolyaBasicProofs}  we give an alternative proof.

According to the experiment described in Section \ref{Models}, the balls that are placed first in the line have an advantage to be picked (the white if $q\in(0, 1)$, the black if $q>1$). In fact, this leads to the extinction of drawings from the balls of the other color; there is a point after which the number of balls in the urn of that color stays fixed to a random number. In the next theorem, we identify the distribution of this number. We treat the case $q>1$.

\begin{theorem}[Extinction of the second color] \label{LimitOfWeakColor} Assume that $q>1, r\in \N\cup\{\infty\}, s\in \N$. As $n\to\infty$,  with probability one, $(X_n)_{n\ge1}$ converges to a random variable $X$ with values in $\N$ and probability mass function

\noindent (i) 
\begin{equation} \label{LimOfWeakColor}
f(x)=
q^{-sx}{\frac{r}{k}+x-1 \brack x}_{q^{-k}}\frac{(q^{-s}; q^{-k})_\infty}{(q^{-r-s}; q^{-k})_\infty}
\end{equation} 
for all $x\in\N$ in the case $r\in \N$ and

\noindent (ii) 
\begin{equation}
f(x)=\left(\frac{q^{-s}}{1-q^{-k}}\right)^x \frac{1}{[x]_{q^{-k}}!}   (q^{-s}; q^{-k})_\infty 
\end{equation}
for all $x\in\N$ in the case $r=\infty$.
\end{theorem}
When $r\in\N$ and $k|r$, $X$ has the negative $q$-binomial distribution of the second kind with parameters $r/k, q^{-s}, q^{-k}$ (see \S 3.1 in \cite{Char16} for its definition). When $r=\infty$, $X$ has the Euler distribution with parameters $q^{-s}/(1-q^{-k}), q^{-k}$ (see \S 3.3 in \cite{Char16} again).

\subsection{$q$-P\'olya urn. Scaling limits}

As in Section \ref{Polya_Scaling_Sec}, we consider an urn whose composition after $n$ drawings is $A_n^{(m)}$ white and $B_n^{(m)}$ black balls. $m\in \Np$ is a parameter.  Our objective is to find limits of the entire path of the process $(A_n^{(m)})_{n\in \N}$ analogous to the ones of  Section \ref{Polya_Scaling_Sec} for the P\'olya urn. Assume that $q>1$.

If we keep $q$ fixed, nothing new appears because: (a) If $A_0^{(m)}, B_0^{(m)}$  are fixed for all $m$, then after some point  we pick only black balls (Theorem \ref{LimitOfWeakColor}(i)). (b) If $\lim_{m\to\infty} B_0^{(m)}=\infty$ then the process converges to the one where we pick only black balls. (c) If  $B_0^{(m)}$ is fixed for all $m$ and $\lim_{m\to\infty} A_0^{(m)}=\infty$ then the process converges to the one where $r=\infty$ and again, after some point, we pick only black balls (Theorem \ref{LimitOfWeakColor}(ii)).

Interesting limits appear once we take $q=q_m$ to depend on $m$ and approach 1 as $m\to\infty$. We study two regimes for $q_m$. In the first, the distance of $q_m$ from 1 is $\Theta(1/m)$ while in the second, the distance is $o(1/m)$.

\subsubsection{The regime $q=1+\Theta(m^{-1})$}

Assume that $q_m=c^{1/m}$ with $c>1$.
\begin{theorem}  \label{QPolyaPathRegime1}  Fix  $w_0\in \Np$ and $b_0\ge 0$. If $A^{(m)}_0=w_0$ and $\lim_{m\to\infty} B^{(m)}_0/m=b_0$, then the process  $(k^{-1}(A^{(m)}_{[m t]}-A_0^{(m)}))_{t\ge0}$ converges in distribution as $m\to\infty$ to an inhomogeneous in time pure birth process $Z$ with starting value $0$ and such that for all $0\le t_1<t_2, j\in\N$, the random variable $Z(t_2)-Z(t_1)|Z(t_1)=j$ has distribution $NB\big(\frac{w_0}{k}+j, \frac{1-c^{-b_0-kt_1}}{1-c^{-b_0-kt_1}}\big)$. Equivalently, $Z$ has rates 
\begin{equation} \label{Qrates}
\gl_{t, j}=\frac{w_0+j k}{c^{b_0+kt}-1}\log c 
\end{equation}
for all $(t, j)\in[0, \infty)\times \N$. 
\end{theorem}

\begin{theorem}  \label{QPolyaPathRegime2} Assume that $A^{(m)}_0=g_m$ and $\lim_{m\to\infty} B^{(m)}_0/m=b_0$,  where $b_0\in(0, \infty)$ and $g_m\in \Np, g_m\to\infty, g_m=o(m)$ as $m\to\infty$. Then the process  $(k^{-1}(A^{(m)}_{[t m/g_m]}-A^{(m)}_0))_{t\ge0}$ converges in distribution, as $m\to\infty$, to the Poisson process on $[0, \infty)$ with rate 
\begin{equation} \label{QPrates}
\frac{\log c}{c^{b_0}-1}.
\end{equation}
\end{theorem}

\begin{theorem} \label{QPolyaPathODE}
Assume that $A_{0}^{(m)}, B_{0}^{(m)}$ are such that $\lim_{m\to\infty} A_{0}^{(m)}/m=a, \lim_{m\to\infty}B_{0}^{(m)}/m=b$, where $a, b\in[0, \infty)$ are not both zero. Then the process $\left(A_{[mt]}/m\right)_t\geq 0$ converges in distribution, as $m\rightarrow+\infty$, to the unique solution of the differential equation 
\begin{align}
\hat X_{0}&=a, \\
d\hat X_{t}&=k\frac{1-c^{\hat X_{t}}}{1-c^{a+b+k t}}dt, \label{ODEqPolya}
\end{align}
which is 
\beq \label{qPolyaDetLimit}
\hat X_t:=a-\frac{1}{\log c}\log\left(\frac{c^b-1+c^{-kt}(1-c^{-a})}{c^b-c^{-a}} \right).
\eeq
\end{theorem}

As for the P\'olya urn, we determine the fluctuations of the process  $(A_{[mt]}^{(m)}/m)_{t\ge0}$ around its $m\to\infty$ limit, $\hat X$. Let 
$$\hat C_t^{(m)}=\sqrt{m}\bigg( \frac{A_{[mt]}^{(m)}}{m}-\hat X_t\bigg)$$
for all $m\in\Np$ and $t\ge0$.

\begin{theorem} \label{QPolyaDiffusionLimit} 
Let $a, b\in[0, \infty)$, not both zero, $\theta_1, \theta_2\in \R$,  and assume that $A_0^{(m)}:=[am+\theta_1\sqrt{m}], B_0^{(m)}=[bm+\theta_2\sqrt{m}]$ for all large $m\in\N$. Then the process $(\hat C_t^{(m)})_{t\ge0}$ converges in distribution, as $m\to\infty$, to the unique solution of
the stochastic differential equation
\begin{equation} \label{QPolyaNoiseSDE}
\begin{aligned}
\hat Y_0&=\theta_1, \\
d\hat Y_t&=\frac{k\log c}{c^{a+b+k t}-1} \bigg\{\frac{(c^{a+b}-1)\hat Y_t-c^b(c^a-1)(\theta_1+\theta_2)}{c^b-1+c^{-kt}(1-c^{-a})}\bigg\} \,dt\\&+k\sqrt{(c^a-1)(c^b-1)} \frac{c^{(a+kt)/2}}{c^{a+b+kt}-c^{a+kt}+c^a-1}\,dW_t,
\end{aligned}
\end{equation}
which is
\begin{equation} \label{QPolyaDiffLimit}
\begin{aligned} 
\hat Y_t=\frac{c^{a+b+kt}-1}{c^{a+b+kt}-c^{a+kt}+c^a-1}\bigg(&\theta_1-(\theta_1+\theta_2)\frac{c^{a+b}(c^a-1)}{c^{a+b}-1}\frac{c^{kt}-1}{c^{a+b+kt}-1}\\&+k\sqrt{(c^a-1)(c^b-1)} \int_0^t\frac{c^{(a+kt)/2}}{c^{a+b+kt}-1}\, dW_s \bigg).
\end{aligned}
\end{equation}
$W$ is a standard Brownian motion
\end{theorem}

\subsubsection{The regime $q=1+o(m^{-1})$}

In this regime, we let $q=q(m):=c^{\gep_m/m}$ where $c>1$ and $\gep_m\to0^+$ as $m\to\infty$. With computations analogous to those of the results of the previous subsection, it is easy to see that Theorems \ref{PolyaPathRegime1}, \ref{PolyaPathRegime2} , \ref{PolyaPathLinear}, \ref{PolyaDiffusionLimit} hold exactly the same for the $q$-P\'olya urn in this regime.

\subsection{$q$-P\'olya urn with many colors}

In this paragraph, we give a $q$-analog for the P\'olya urn with more than two colors. The way to do the generalization is inspired by the experiment we used in order to explain relation \eqref{qPWhite}.

Let $l\in\N, l\ge 2$, and $q\in(0,1)$. Assume that we have an urn containing $w_i$ balls of color $i$ for each $i\in\{1, 2, \ldots, l\}$. To draw  a ball from the urn, we do the following. We order the balls in a line, first those of color 1, then those of color 2, and so on. Then we visit the balls, one after the other, in the order that they have been placed, and we select each with probability $1-q$ independently of what happened with the previous balls. If we go through all balls without picking any, we repeat the same procedure starting from the beginning of the line. Once a ball is selected, the drawing is completed. We return the ball to its position together with another $k$ of the same color. For each $i=0, 1, \ldots, l$, let $s_i=\sum_{1\le j\le i} w_j$. Notice that $s_l$ is the total number of balls in the urn. Then, working as for \eqref{LinePick}, we see that
\begin{equation} \label{problcolors}
\PP(\text{color $i$ is drawn})=q^{s_{i-1}}\frac{1-q^{w_i}}{1-q^{s_l}}=\frac{q^{s_{i-1}}-q^{s_i}}{1-q^{s_l}}=q^{s_{i-1}}\frac{[w_i]_q}{[s_l]_q}. 
\end{equation} 

Call $p_i$ the number in the last display for all $i=1, 2, \ldots, l$. Note that when $q\to1$, $p_i$ converges to $w_i/s_l$, which is the probability for the usual P\'olya urn with $l$ colors.  It is clear that for any given $q\in(0, \infty)\sm\{1\}$, the numbers $p_1, p_2, \ldots, p_l$ are non-negative and add to 1 (the second fraction in \eqref{problcolors} shows this). We define then for this $q$ the $q$-P\'olya urn with colors $1, 2, \ldots, l$ to be the sequential procedure in which, at each step, we add $k$ balls of a color picked randomly among $\{1, 2, \ldots, l\}$ so that the probability that this color is $i$ is $p_i$ . 

When $q>1$, these probabilities come out of the experiment described above but in which we place the balls in reverse order  (that is, first those of color $l$, then those of color $l-1$, and so on) and we go through the list selecting each ball with probability $1-q^{-1}$. It is then easy to see that the probability to pick a ball of color $i$ is $p_i$.

\begin{theorem} \label{PMFOfWeakColors} Assume that $q\in(0, 1)$ and that we start with $a_1, a_2,\ldots, a_l$ balls from colors 1, 2, $\ldots, l$ respectively, where $a_1, a_2, \ldots, a_l\in\N$ are not all zero. Call $X_{n, i}$ the number of times in the first $n$ drawings that we picked color $i$. The probability mass function for the vector $\left(X_{n,2},X_{n,3},\ldots,X_{n,l}\right)$ is 
\begin{align} \label{probmassforvectorAig}
\PP\left(X_{n,2}=x_{2},\ldots,X_{n,l}=x_{l}\right)&=q^{\sum_{i=2}^{l}x_i\sum_{j=1}^{i-1}\left(a_{j}+k x_{j}\right)}\frac{\prod_{i=1}^{l}{-\frac{a_{i}}{k} \brack x_{i}}_{q^{-k}}}{{-\frac{a_1+a_2\ldots +a_l}{k} \brack n}_{q^{-k}}}
\\  \label{probmassforvector} 
&={n \brack x_{1},x_{2},\ldots,x_{l}}_{q^{-k}}  \frac{q^{\sum_{i=2}^{l}x_i\sum_{j=1}^{i-1}\left(a_{j}+k x_{j}\right)}\prod_{i=1}^{l}\left[-\frac{a_{i}}{k}\right]_{x_{i},q^{-k}}}{\left[-\frac{a_{1}+a_{2}+\ldots+a_{l}}{k}\right]_{n,q^{-k}}}
\end{align}
for all $x_2, \ldots, x_l\in\{0,1,2,\ldots,n\}$ with $x_2+\cdot+x_l\le n$, where $x_{1}:=n-\sum_{i=2}^{l}x_{2}$ and ${n \brack x_{1},x_{2},\ldots,x_{l}}_{q^{-k}}:=\frac{[n]_{q^{-k}}!}{[x_{1}]_{q^{-k}}!\cdot \ldots \cdot [x_{l}]_{q^{-k}}!}$ is the $q$-multinomial coefficient.
\end{theorem}

It follows from Theorem \ref{LimitOfWeakColor} that when $q\in(0, 1)$, after some random time, we will be picking only balls of color 1. So that the number of times that we pick each of the other colors $2, 3, \ldots, l$, say $X_2, X_3, \ldots, X_n$ are finite. We determine the joint distribution of these numbers.

\begin{theorem}\label{LimitOfWeakColors} Under the assumptions of Theorem \ref{PMFOfWeakColors},
as $n\rightarrow +\infty,$ with probability one, the vector $\left(X_{n,2},X_{n,3},\ldots,X_{n,l}\right)$ converges  to a random vector $\left(X_{2},X_{3},\ldots,X_{l}\right)$ with values in $\mathbb{N}^{l-1}$ and probability mass function 
\begin{equation}
f\left(x_{2},x_{3},\ldots,x_{l}\right)=q^{\sum_{i=2}^{l}x_{i}\sum_{j=1}^{i-1}a_{j}}\prod_{i=2}^{l}{x_{i}+\frac{a_{i}}{k}-1 \brack x_{i}}_{q^{k}} \frac{(q^{a_1}; q^k)_\infty}{(q^{a_1+\cdots+a_l}; q^k)_\infty}
\end{equation}
for all $x_2, \ldots, x_l\in \N$.
\end{theorem}
\noindent Note that the random variables $X_2, \ldots, X_l$ are independent although $\left(X_{n,2},X_{n,3},\ldots,X_{n,l}\right)$ are dependent.

Next, we look for a scaling limit for the path of the process. Assume that $c\in(0, 1)$ and $q_{m}=c^{1/m}$.  Let $A_{j,i}^{(m)}$ be the number of balls of color $i$ in this urn after $j$ drawings.
\begin{theorem} \label{QPolyaPathODEWColors}
Let $m$ be a positive integer and assume that in the $q$-P\'olya urn with $l$ different colors of balls it holds  $\frac{1}{m}\left(A_{0,1}^{(m)},A_{0,2}^{(m)}, \ldots , A_{0,l}^{(m)}\right)\overset{m\to\infty}{\to}\left(a_{1},  a_{2} ,\dots, a_{l}\right)$, where $a_{1},\ldots,a_{l}\in\left[0,\infty\right)$ are not all zero. Set $\sigma_0=0$ and $\sigma_i:=\sum_{j\le i}a_j$ for all $i=1, 2, \ldots, l$. Then the process $\left(\frac{1}{m}A_{[mt],1}^{(m)},\frac{1}{m}A_{[mt],2}^{(m)},\ldots,\frac{1}{m}A_{[mt],l}^{(m)}\right)_{t\geq 0} $ converges in distribution, as $m\rightarrow +\infty$, to $(X_{t, 1}, X_{t, 2}, \ldots, X_{t, l})_{t\ge0}$ with
\begin{equation} \label{SystemSolution}
X_{t, i}=a_i+\frac{1}{\log c} \log \frac{(1-c^{\sigma_l+kt})-c^{\sigma_{i-1}}(1-c^{k t})}{(1-c^{\sigma_l+kt})-c^{\sigma_i}(1-c^{k t})}
\end{equation}
for all $i=1, 2, \ldots, l$.
\end{theorem}

As in the case of two colors, we study the regime where $q_{m}=c^{\epsilon_{m}/m}$, with $c\in(0,1)$ and $\epsilon_{m}\rightarrow 0^{+}.$
\begin{theorem} \label{QPolyaPathODEWColors2}
Let $m$ be a positive integer and assume that in the q-P\'olya urn with $l$ different colors of balls that  $\frac{1}{m}\left(A_{0,1}^{(m)},A_{0,2}^{(m)}, \ldots , A_{0,l}^{(m)}\right)\overset{m\to\infty}{\to}\left(a_{1},  a_{2} ,\dots, a_{l}\right)$, where $a_{1},\ldots,a_{l}\in\left[0,\infty\right)$ are not all zero.  Then the process $\left(\frac{1}{m}A_{[mt],1}^{(m)},\frac{1}{m}A_{[mt],2}^{(m)},\ldots,\frac{1}{m}A_{[mt],l}^{(m)}\right)_{t\geq 0}$ converges in distribution, as $m\rightarrow +\infty$, to $(X_t)_{t\ge0}$ with
\begin{equation}\label{SystemSolutionReg2}  X_t=\left(1+\frac{kt}{a_{1}+\dots+a_{l}}\right)(a_1, a_2, \ldots, a_l)
\end{equation}
for all $t\ge0$.
\end{theorem}

\noindent \textbf{Orientation}. In Section \ref{qPolyaBasicProofs}, we prove Fact 1 and Theorem 1.5, which are basic results for the $q$-P\'olya urn. Section \ref{JumpLimits} (Section \ref{DetermDiffLimits}) contains the proofs of the theorems for the P\'olya and $q$-P\'olya urns that give convergence to a jump process (to a continuous process). Finally, Section \ref{ProofMoreColors} contains the proofs for the results that refer to the $q$-P\'olya urn with arbitrary, finite number of colors.

\section{$q$-P\'olya urn. Prevalence of a single color} \label{qPolyaBasicProofs}

In this section, we prove the claims of Section \ref{qPolyaBasic}. Before doing so, we mention three properties of the $q$-binomial coefficient.  For all $q\in (0, \infty)\sm\{1\}, x\in \mathbb{C}, n, k\in \N$ with $k\le n$ it holds 
\begin{align}
\label{NegToPosa}
[-x]_q&=-q^{-x}[x]_q,\\
\label{NegToPos}
{-x \brack k}_q&=(-1)^k q^{-k(k+2x-1)/2}{x+k-1 \brack k}_q,\\ 
{x \brack k}_{q^{-1}}&=q^{-k(x-k)}{x \brack k}_{q}   \label{qMinusOneToq},\\
\sum_{1\le i_1<i_2<\cdots<i_k\le n} q^{i_1+i_2+\cdots+i_k}&=q^{{k+1 \choose 2}}{n \brack k }_{q}. \label{QBinCoeffCount}
\end{align}
The first is trivial, the second follows from the first, the third is easily shown, while the last is Theorem 6.1 in \cite{KacCheung}.

\begin{proof}[\textbf{Proof of Fact 1}] (i) The probability to get black balls exactly at the drawings $i_1<i_2<\cdots<i_{n-x}$ is
\begin{equation} \label{PnWBalls}  g(i_1, i_2, \ldots, i_{n-x})=\frac{\prod_{j=0}^{x-1} [r+j k]_q \prod_{j=0}^{n-x-1} [s+j k]_q}{\prod_{j=0}^{n-1} [r+s+j k]_q} q^{\sum_{\nu=1}^{n-x}r+(i_\nu-\nu) k}.
\end{equation}
To see this, note that, due to \eqref{qPWhite} and \eqref{qPBlack}, the required probability would be equal to the above fraction if in \eqref{qPBlack} the term $q^w$ were absent. This term appears whenever we draw a black ball. Now, when we draw the $\nu$-th black ball, there are $r+(i_\nu-\nu) k$ white balls in the urn, and this explains the exponent of $q$ in \eqref{PnWBalls}.

Since $[x+j k]_q=\frac{1-q^{x+jk}}{1-q}=[-\frac{x}{k}-j]_{q^{-k}}[-k]_q$ for all $x, j\in \R$, the fraction in \eqref{PnWBalls} equals 
\begin{equation}
\frac{[-a]_{x, q^{-k}} [-b]_{n-x, q^{-k}}}{[-a-b]_{n , q^{-k}}}.
\end{equation}
Then
\begin{align}\sum_{1\le i_1<i_2<\cdots<i_{n-x}\le n} q^{\sum_{\nu=1}^{n-x}r+(i_\nu-\nu)k}&=q^{(n-x) r-k(n-x)(n-x+1)/2}\sum_{1\le i_1<i_2<\cdots<i_{n-x}\le n} (q^k)^{i_1+i_2+\cdots+i_{n-x}}\\
&=q^{(n-x) r-k(n-x)(n-x+1)/2} q^{k {n-x+1 \choose 2 }} {n \brack x }_{q^k} \\
&=q^{(n-x) r} q^{kx(n-x)}{n \brack x }_{q^{-k}}=q^{k (n-x)(a+x)}{n \brack x }_{q^{-k}}.
\end{align}
The second equality follows from \eqref{QBinCoeffCount} and the equality  ${n \brack x }_{q^k}={n \brack n-x }_{q^k}$. The third, from \eqref{qMinusOneToq}. 
Thus, the sum $\sum_{1\le i_1<i_2<\cdots<i_{n-x}\le n}g(i_1, i_2, \ldots, i_{n-x})$ equals the first expression in \eqref{QPolyaWhiteDistr}. The second expression in \eqref{QPolyaWhiteDistr} and \eqref{QPWDQ} follow by using \eqref{NegToPos} and \eqref{qMinusOneToq} respectively. 

(ii)  In this scenario, we take $r\to\infty$ in the last expression in \eqref{QPolyaWhiteDistr}. We will explain shortly why this gives the probability we want. Since $q^{-k}\in(0, 1)$, we have $\lim_{t\to\infty}[t]_{q^{-k}}=(1-q^{-k})^{-1}$  and thus, for each $\nu\in\N$, it holds 
\begin{equation}
\lim_{t\to\infty}{t+\nu-1 \brack \nu}_{q^{-k}}=\frac{1}{[\nu]_{q^{-k}}!} \frac{1}{(1-q^{-k})^\nu}.
\end{equation}
Applying this twice in the last expression in \eqref{QPolyaWhiteDistr} (there $a=r/k\to\infty$), we get as limit the right hand side of \eqref{QPolyaWhiteDistrRInfty}.

Now, to justify that passage to the limit $r\to\infty$ in \eqref{QPolyaWhiteDistr} gives the required result, we argue as follows.
For clarity, denote the probability $\PP_q(\text{white})$ when there are $w$ white and $b$ black balls in the urn by $\PP_q^{w, b}(\text{white})$. And when there are $r$ white and $s$ black balls in the urn in the beginning of the procedure, denote the probability of the event $X_n=x$ by $\PP^{r, s}(X_n=x)$. It is clear that the probability $\PP^{r, s}(X_n=x)$ is a continuous function (in fact, a polynomial) of the quantities
$$\PP_q^{r+ki, s+k j}(\text{white}):i=0, 1, \ldots, x-1, j=0, 1, \ldots, n-x-1,$$
for all values of $r\in \N\cup\{\infty\}, s\in \N$. In $\PP^{\infty, s}(X_n=x)$, each such quantity, $\PP_q^{\infty, m}(\text{white})$, equals $\lim_{r\to\infty}\PP^{r, m}(\text{white})$. Thus,  $\PP^{\infty, s}(X_n=x)=\lim_{r\to\infty}\PP^{r, s}(X_n=x)$.  
\end{proof}

Before proving Theorem \ref{LimitOfWeakColor}, we give a simple argument that shows that eventually we will be picking only black balls. That is, the number $X:=\lim_{n\to\infty} X_n$ of white balls drawn in an infinite sequence of drawings is finite. It is enough to show it in the case that $r=\infty$ and $s=1$ since, by the experiment that realizes the $q$-P\'olya urn, we have (using the notation from the proof of Fact 1 (ii))
$$\PP^{r, s}(X=\infty)\le \PP^{\infty, 1}(X=\infty).$$
For each $n\in\Np$, call $E_n$ the event that at the $n$-th drawing we pick a white ball, $B_n$ the number of black balls present in the urn after that drawing (also, $B_0:=1$), and write $\hat q:=1/q$. Then $\PP(E_n)=\EE(\PP(E_n| B_{n-1}))=\EE(\hat q^{B_{n-1}})$. We will show that this decays exponentially with $n$. Indeed, since at every drawing there is probability at least $1-\hat q$ to pick a black ball, we can construct in the same probability space the random variables $(B_n)_{n\ge1}$ and $(Y_i)_{i\ge1}$ so that the $Y_i$ are i.i.d. with $Y_1\sim$ Bernoulli$(1-\hat q)$ and $B_n\ge 1+k(Y_1+\cdots+Y_n)$ for all $n\in\Np$. Consequently, 
$$\PP(E_n)\le \EE(\hat q^{1+k(Y_1+\cdots+X_{n-1})})=\hat q \{\EE(\hat q^{kY_1})\}^{n-1}.$$
This implies that $\sum_{n=1}^\infty \PP(E_n)<\infty$, and the first Borel-Cantelli lemma gives that  $\PP^{\infty, 1}(X_\infty=\infty)=0$.

\begin{proof}[\textbf{Proof of Theorem \ref{LimitOfWeakColor}}]  Since $(X_n)_{n\ge1}$ is increasing, it converges to a random variable $X$ with values in $\N\cup\{\infty\}$. In particular, it converges to this variable in distribution. 
Our aim is to take the limit as  $n\to\infty$ in the last expression in \eqref{QPolyaWhiteDistr} and in \eqref{QPolyaWhiteDistrRInfty} in order to determine the distribution of $X$. Note that for $a\in \R$  and $\theta\in[0, 1)$ it is immediate that (recall \eqref{qPoch} for the notation)
\begin{equation} \label{QBinomialAsymptotics}
\lim_{n\to\infty}{a+n \brack n }_{\theta}=\frac{(\theta^{a+1}; \theta)_\infty}{(\theta; \theta)_\infty}.
\end{equation}

\noindent (i) Taking $n\to\infty$ in the last expression in \eqref{QPolyaWhiteDistr} and using \eqref{QBinomialAsymptotics}, we get the required expression, \eqref{LimOfWeakColor}, for $f$. Then relation (2.2) in \cite{Char12} (or (8.1) in \cite{KacCheung}) shows that $\sum_{x\in\N} f(x)=1$, so that it is a probability mass function of a random variable $X$ with values in $\N$.

\noindent (ii) This follows after taking limit in \eqref{QPolyaWhiteDistrRInfty} and using \eqref{QBinomialAsymptotics} and $\lim_{n\to\infty} (1-q^{-k})^n [n]_{q^{-k}}!=(q^{-k}; q^{-k})_\infty$.
\end{proof}

\section{Jump process limits. Proof of Theorems \ref{PolyaPathRegime1}, \ref{PolyaPathRegime2}, \ref{QPolyaPathRegime1}, \ref{QPolyaPathRegime2}} \label{JumpLimits}

In the case of Theorems \ref{PolyaPathRegime1}, \ref{QPolyaPathRegime1}, we let $g_m:=1$ for all $m\in \N^+$, and in all four theorems we let $v:=v_m:=m/g_m$. Our interest is in the sequence of the processes $(Z^{(m)})_{m\ge1}$ with
\begin{equation}
Z^{(m)}(t)=\frac{1}{k}(A^{(m)}_{[v t]}-A^{(m)}_0)
\end{equation}
for all $t\ge0$.

We apply Theorem 7.8 in \cite{EtKu}, that is, we show that the sequence $(Z^{(m)})_{m\ge1}$ is tight and its finite dimensional distributions converge. Tightness gives that there is a subsequence of this sequence that converges in distribution to a process $Z=(Z_t)_{t\ge0}$ with paths in the space $D_\R[0, \infty)$ of real valued functions on $[0, \infty)$ right continuous with left limits. Then tightness together with convergence of finite dimensional distributions shows that the whole sequence $(Z^{(m)})_{m\ge1}$ converges in distribution to $Z$.

\medskip

\noindent \textbf{Notation:} (i) For sequences $(a_n)_{n\in\N}, (b_n)_{n\in\N}$ with values in $\R$, we will say that they are asymptotically equivalent, and will write $a_n\sim b_n$ as $n\to\infty$, if $\lim_{n\to\infty} a_n/b_n=1$. We use the same expressions for functions $f, g$ defined in a neighborhood of $\infty$ and satisfy $\lim_{x\to\infty} f(x)/g(x)=1$.

\noindent (ii) For $a\in \mathbb{C}$ and $k\in \Np$, let
\begin{align}
(a)_k&:=a(a-1)\cdots (a-k+1),\\
a^{(k)}&:=a(a+1)\cdots(a+k-1),
\end{align}
the falling and rising factorial respectively. Also let $(a)_0:=a^{(0)}:=1$.

\subsection{Convergence of finite dimensional distributions} \label{FddConvergence}

Since for each $m\ge1$ the process $Z^{(m)}$ is Markov taking values in $\N$ and increasing in time, it is enough to show that the conditional probability  
\begin{equation}\label{ConditionalTransition}  \PP(Z^{(m)}(t_2)=k_2 | Z^{(m)}(t_1)=k_1)
\end{equation}
converges as $m\to\infty$ for each $0\le t_1<t_2$ and nonnegative integers $k_1\le k_2$.

\noindent Consider first the case of P\'olya urn and define 
\begin{align} 
n&:=[v t_2]-[v t_1], \label{nDef}\\
x&:=k_2-k_1,\\
\sigma&:=\frac{A_0^{(m)}+kk_1}{k},\\
\tau&:=\frac{k[v t_1]-kk_1+B_0^{(m)}}{k}. \label{tauDef}
\end{align}
Then, the above probability equals
\begin{align} \notag
&\PP(A_{[v t_2]}^{(m)}=kk_2+w_0 | A_{[v t_1]}^{(m)}=kk_1+w_0)\\
&={n \choose x}\frac{k\sigma (k\sigma+k)(k\sigma+2k)\cdots (k\sigma+(x-1)k) k\tau(k\tau+k)(k\tau+2k)\cdots (k\tau+(n-x-1)k)}{(k\sigma+k\tau)(k\sigma+k\tau+k)(k\sigma+k\tau+2k)\cdots (k\sigma+k\tau+(n-1)k)}\\
&=\frac{(n)_x}{x!}\frac{\sigma^{(x)} \tau^{(n-x)} }{(\sigma+\tau)^{(n)}}=\frac{(n)_x}{x!} \sigma^{(x)} \frac{\Gamma\left(\tau+n-x\right)}{\Gamma\left(\tau\right)}\frac{\Gamma\left(\sigma+\tau\right)}{\Gamma\left(\sigma+\tau+n\right)}. \label{TransitionProbab}
\end{align}
To compute the limit as $m\to\infty$ of \eqref{TransitionProbab}, we will use Stirling's approximation for the Gamma function,
\begin{equation} \label{GammaStirling}
\Gamma(y)\sim\left(\frac{y}{e}\right)^y\sqrt{\frac{2\pi}{y}}
\end{equation}
as $y\to\infty$, and its consequence
\begin{equation} \label{GammaAsymptotics}
\Gamma(y+a)\sim\Gamma(y) y^a
\end{equation}
as $y\to\infty$ for all $a\in\mathbb{R}$. 

\begin{proof}[\textbf{Theorem \ref{PolyaPathRegime1}}] Recall that $v=m$ in this theorem. Using \eqref{GammaAsymptotics} twice, with the role of $a$ played by $-x$ and $\sigma$, we see that the last quantity in \eqref{TransitionProbab}, for $m\to\infty$, is asymptotically equivalent to 
\begin{align*}
&\frac{(m(t_2-t_1))^x}{x!}\sigma^{(x)}\tau^\sigma \frac{\left(\tau+n\right)^{-x}}{\left(\tau+n\right)^{\sigma}}\sim \frac{(m(t_2-t_1))^x}{x!}\sigma^{(x)} \frac{\{m(t_1+(b_0/k))\}^\sigma}{\{m(t_2+(b_0/k))\}^{\sigma+x}} \\
&=\frac{(t_2-t_1)^x}{x!} \sigma^{(x)}  \frac{\{t_1+(b_0/k)\}^\sigma}{\{t_2+(b_0/k)\}^{\sigma+x}} ={\sigma+x-1 \choose x} \left(\frac{t_2-t_1}{t_2+(b_0/k)}\right)^x \left(1-\frac{t_2-t_1}{t_2+(b_0/k)}\right)^\sigma.
\end{align*}
Thus, as $m\to\infty$, the distribution of $\{Z^{(m)}(t_2)-Z^{(m)}(t_1)\}| Z^{(m)}(t_1)=k_1$ converges to the negative binomial distribution with parameters $\sigma, \frac{t_1+(b_0/k)}{t_2+(b_0/k)}$ [recall \eqref{NBDensity}]. 
\end{proof}

\begin{proof}[\textbf{Theorem \ref{PolyaPathRegime2}}]
Using \eqref{GammaStirling}, we see that the last quantity in \eqref{TransitionProbab}, for $m\to\infty$, is asymptotically equivalent to 
\begin{align*}
&\frac{(m(t_2-t_1))^x}{x! g_m^x}\frac{g_m^x}{k^x} e^x \frac{(\tau+n-x)^{\tau+n-x}}{\tau^\tau}\frac{(\sigma+\tau)^{\sigma+\tau}}{(\sigma+\tau+n)^{\sigma+\tau+n}}
 \\
&\sim  \frac{m^x (t_2-t_1)^x}{x! k^x} e^x (\tau+n-x)^{-x} \left(\frac{\tau+n-x}{\sigma+\tau+n}\right)^n  \left(\frac{\sigma+\tau}{\sigma+\tau+n}\right)^\sigma   \left(\frac{(\tau+n-x)(\sigma+\tau)}{\tau(\sigma+\tau+n)}\right)^{\tau} \\
&\sim  \frac{m^x (t_2-t_1)^x}{x! k^x} e^x  \tau^{-x} e^{-(t_2-t_1)/b_0}   e^{-(t_2-t_1)/b_0}  e^{-x+(t_2-t_1)/b_0}\sim  \frac{1}{x!}\left(\frac{t_2-t_1}{b_0}\right)^x   e^{-(t_2-t_1)/b_0} . 
\end{align*}
Here it was crucial that $b_0>0$.
Thus, as $m\to\infty$, the distribution of $\{Z^{(m)}(t_2)-Z^{(m)}(t_1)\}| Z^{(m)}(t_1)=k_1$ converges to the Poisson distribution with parameter $(t_2-t_1)/b_0$.
\end{proof}

\noindent Now we treat the cases of Theorems \ref{QPolyaPathRegime1}, \ref{QPolyaPathRegime2}, which concern the $q$-P\'olya urn.
Define again $n, x, \sigma, \tau$ as in \eqref{nDef}-\eqref{tauDef}, and $r:=q_m^{-k}=c^{-k/m}$.
Then, the probability in \eqref{ConditionalTransition}, with the help of the last expression in \eqref{QPolyaWhiteDistr}, is computed as
\beq \label{QPolyaTransition}
r^{\tau x} {\sigma+x-1 \brack x }_{r}  
\frac{{\tau+n-x-1 \brack n-x }_{r} }{{\sigma+\tau+n-1 \brack n }_{r} }=r^{\tau x} {\sigma+x-1 \brack x }_{r}  \Big(\prod_{i=n-x+1}^{n}(1-r^i)\Big)   \frac{1}{\prod_{i=n-x}^{n-1}(1-r^{\tau+i})} \frac{[\tau+n-1]_{n, r}}{[\sigma+\tau+n-1]_{n, r}}.
\eeq 
The last ratio is 
\begin{align}
\prod_{i=0}^{n-1} \frac{1-r^{\tau+i}}{1-r^{\sigma+\tau+i}}=\prod_{i=0}^{n-1} \left(1-(1-r^\sigma) r^\tau \frac{r^i}{1-r^{\sigma+\tau+i}}\right).
\end{align}
Denote by $1-a_{m, i}$ the $i$-th term of the product. The logarithm of the product equals 
\begin{equation} \label{RiemannSum}
-(1-r^\sigma) r^\tau  \sum_{i=0}^{n-1} \frac{r^i}{1-r^{\sigma+\tau+i}}+o(1)
\end{equation}
as $m\to\infty$. To justify this, note that $1-r^\sigma\sim \frac{1}{m}(A_0^{(m)}+k k_1) \log c$ and $r^{\tau+i}/(1-r^{\sigma+\tau+i})\le 1/(1-c^{-b_0})$ for all $i\in \N$. Thus, for all large $m$, $|a_{m, i}|<1/2$ for all $i=0, 1, \ldots, n-1$, and the error in approximating the logarithm of $1-a_{m, i}$ by $-a_{m, i}$ is at most $|a_{m, i}|^2$ (by Taylor's expansion, we have $|\log(1-y)+y|\le |y|^2$ for all $y$ with $|y|\le1/2$). The sum of all errors is at most $n \max_{0\le i<n}|a_{m, i}|^2$, which goes to zero as $m\to\infty$ because $1-r^\sigma\sim C/n$ for some appropriate constant $C>0$.

We will compute the limit of the right hand side of \eqref{QPolyaTransition} as $m\to\infty$ under the assumptions of Theorems \ref{QPolyaPathRegime1}, \ref{QPolyaPathRegime2}.

\begin{proof}[\textbf{Theorem \ref{QPolyaPathRegime1}}]
As $m\to\infty$, the first term of the product in \eqref{QPolyaTransition} converges to $c^{-x(b_0+k t_1)}$. The $q$-binomial coefficient converges to ${k^{-1}w_0+k_2-1 \choose k_2-k_1}$. The third term converges to $(1-c^{-k(t_2-t_1)})^x$,  while the denominator of the fourth term converges to $(1-\rho_2)^x$, where we set $\rho_i:=c^{-b_0-kt_i}$ for $i=1, 2$. The expression preceding $o(1)$ in \eqref{RiemannSum} is asymptotically equivalent to
\begin{align}&-\frac{k}{m} \sigma (\log c) \rho_1 \sum_{i=0}^{n-1} \frac{c^{-ki/m}}{1-r^{\sigma+\tau} c^{-ki/m}}
= -\rho_1k \sigma (\log c) \frac{1}{m}\sum_{i=0}^{n-1} \frac{c^{-ki/m}}{1-\rho_1 c^{-ki/m}}+o(1)
\\&=-\rho_1 k \sigma \log c \int_0^{t_2-t_1} \frac{1}{c^{k y}-\rho_1}\, dy+o(1)
=\sigma \log \frac{1-\rho_1}{1-\rho_2}+o(1).
\end{align}
The equality in the first line is true because $\lim_{m\to\infty}r^{\sigma+\tau}=\rho_1$ and the function $x\mapsto c^{-k i/m}/(1-x c^{-ki/m})$ has derivative bounded uniformly in $i, m$ when $x$ is confined to a compact subset of $[0, 1)$. 
Thus, the limit of \eqref{QPolyaTransition}, as $m\to\infty$, is
\beq 
{\sigma+x-1 \choose x} \left(\frac{\rho_1-\rho_2}{1-\rho_2}\right)^x \left(\frac{1-\rho_1}{1-\rho_2}\right)^\sigma,
\eeq
which means that, as $m\to\infty$, the distribution of $\{Z^{(m)}(t_2)-Z^{(m)}(t_1)\}| Z^{(m)}(t_1)=k_1$ converges to the negative binomial distribution with parameters $\sigma, (1-\rho_1)/(1-\rho_2)$.
\end{proof}

\begin{proof}[\textbf{Theorem \ref{QPolyaPathRegime2}}] Now the term $r^{\tau x}$ converges to $c^{-x b_0}$, while 
\begin{align} &{\sigma+x-1 \brack x }_{r}\Big(\prod_{i=n-x+1}^{n}(1-r^i)\Big)=\frac{\prod_{i=0}^{x-1}(1-r^{\sigma+i})}{\prod_{i=1}^x (1-r^i)} \Big(\prod_{i=n-x+1}^{n}(1-r^i)\Big)\\ &\sim \frac{\prod_{i=0}^{x-1}(\sigma+i)}{\prod_{i=1}^x i} \frac{((t_2-t_1)k\log c)^x}{g_m^x}\sim\frac{1}{x!}((t_2-t_1)\log c)^x.
\end{align}
The denominator of the fourth term in \eqref{QPolyaTransition}  converges to $(1-c^{-b_0})^x$. The expression in \eqref{RiemannSum} is asymptotically equivalent to
\beq -r^\tau(1-r^\sigma)\sum_{i=0}^{n-1} \frac{r^i}{1-r^{\sigma+\tau+i}}\sim -c^{-b_0} \frac{g_m}{m}\log c \frac{n}{1-c^{-b_0}} \sim -\frac{\log c}{c^{b_0}-1}(t_2-t_1).
\eeq
In the first $\sim$, we used the fact that the terms of the sum, as $m\to\infty$, converge uniformly in $i$ to $(1-c^{-b_0})^{-1}$. Thus, the limit of \eqref{QPolyaTransition}, as $m\to\infty$, is
\beq 
\frac{1}{x!} \left(\frac{\log c}{c^{b_0}-1}(t_2-t_1)\right)^x e^{-\frac{\log c}{c^{b_0}-1}(t_2-t_1) },
\eeq
which means that, as $m\to\infty$, the distribution of $\{Z^{(m)}(t_2)-Z^{(m)}(t_1)\}| Z^{(m)}(t_1)=k_1$ converges to the Poisson distribution with parameter $\frac{t_2-t_1}{c^{b_0}-1}\log c$.
\end{proof}

For use in the following section, we define
\begin{equation} \label{transitionProb}
U(t_1, t_2, k_1, x):=\lim_{m\to\infty}\PP(Z^{(m)}(t_2)=k_1+x | Z^{(m)}(t_1)=k_1)
\end{equation}
for all $0\le t_1\le t_2, k_1\in \N, x\in \N$. The results of this section show that $U$ as a function of $x\in \N$ is a probability mass function of an appropriate random variable with values in $\N$.

\subsection{Tightness} We apply Corollary 7.4 of Chapter 3 in \cite{EtKu}. According to it, it is enough to show that 
\begin{itemize}
\item[(i)] For each $t\ge0$, it holds  $\lim_{R\to\infty} \varlimsup_{m\to\infty}\PP(|Z^{(m)}(t)|\ge R)=0$.

\item[(ii)] For each $T, \gep>0$, it holds $\lim_{\delta\to0} \varlimsup_{m\to\infty} \PP(w'(Z^{(m)}, \delta, T)\ge \gep)=0.$
\end{itemize}
Here, for any function $f:[0, \infty)\to\R$, we define
$$w'(f, \delta, T):=\inf_{\{t_i\}}\max_{i}\sup_{s, t\in[t_{i-1}, t_i)}|f(s)-f(t)|,$$
where the infimum is over all partitions of the form $0=t_0<t_1<\cdots t_{n-1}<T\le t_n$ with $t_i-t_{i-1}>\delta$ for all $i=1, 2, \ldots, n$.

The first requirement holds because $Z^{(m)}(t)$ converges in distribution as we showed in the previous subsection. The second requirement, since $Z^{(m)}$ is a jump process with jump sizes only 1, is equivalent to 
\beq \label{TightnessJumps}
\lim_{\delta\to0^+} \varlimsup_{m\to\infty} \PP(\text{There are at least two jump times of $Z^{(m)}$ in $[0, T]$ with distance}\le \delta)=0.
\eeq
Call $A_{m, \delta}$ the event inside the probability and for $j=1, 2, \ldots, [T/\delta]$ define $I_j:=((j-1)\delta, (j+1)\delta]$.
Then, for each $\ell\in \N$, the probability $\PP(A_{m, \delta}\cap \{Z^{(m)}(T)\le \ell\})$ is bounded above by 
\begin{align}&\sum_{j=1}^{[T/\delta]}\PP( \{Z^{(m)}(T)\le \ell\}\cap\{\text{There are at least two jump times of $Z^{(m)}$ in } I_j\}) \\
&\le  \sum_{j=1}^{[T/\delta]}\PP( \{Z^{(m)}(T)\le \ell\}\cap\{Z^{(m)}((j+1)\delta)-Z^{(m)}((j-1)\delta)\ge 2\})\\
&\le \sum_{j=1}^{[T/\delta]} \max_{0\le \mu \le \ell} \PP(Z^{(m)}((j+1)\delta)-Z^{(m)}((j-1)\delta)\ge 2 | Z^{(m)}((j-1)\delta)=\mu).
\end{align}
The limit of the last quantity as $m\to\infty$, with the use of the function $U$ of \eqref{transitionProb}, is written as
\beq \label{JumpMaxBound}
 \sum_{j=1}^{[T/\delta]} \max_{0\le \mu \le \ell} \sum_{x=2}^\infty U((j-1)\delta, (j+1)\delta, \mu, x) \le\frac{T}{\delta}  \max_{\substack{0\le \mu\le \ell\\ 1\le j\le [T/\delta]}} \sum_{x=2}^\infty U((j-1)\delta, (j+1)\delta, \mu, x).
\eeq

\textsc{Claim:} The max in \eqref{JumpMaxBound} is bounded above by $\delta^2 C(\ell, T)$ for an appropriate constant $C(\ell, T)\in (0, \infty)$ that does not depend on $m$ or $\delta$.

Assuming the claim and taking $m\to\infty$ in $\PP(A_{m, \delta})=\PP(A_{m, \delta}\cap \{Z^{(m)}(T)\le \ell\})+\PP(A_{m, \delta}\cap \{Z^{(m)}(T)>\ell\})$, we get 
$$\varlimsup_{m\to\infty} \PP(A_{m, \delta})\le \delta C(\ell, T)+\varlimsup_{m\to\infty}\PP(\{Z^{(m)}(T)>\ell\}).$$
Now let $\gep>0$. Because of the validity of (i) in the tightness requirements, there is $\ell$ large enough so that the second term is $<\gep$. Fixing this $\ell$ and taking $\delta\to0$ in the inequality, we get \eqref{TightnessJumps}.

\textsc{Proof of the claim:} We establish the above claim for each of  the Theorems \ref{PolyaPathRegime1}, \ref{PolyaPathRegime2}, \ref{QPolyaPathRegime1}, \ref{QPolyaPathRegime2}. We use the following bounds. If $X, Y$ are random variables with $X\sim$ Poisson($\lambda$) and $Y\sim NB(\nu, p)$ then
\begin{align}
\PP(X\ge 2)&\le \lambda^2, \label{PoissonTail}\\
\PP(Y\ge 2)&\le \frac{\nu (\nu+1)}{2} (1-p)^2. \label{NBTail}
\end{align}
The first inequality is elementary, while the second is true because the difference of the two sides $$\PP(Y\ge 2)-\frac{\nu(\nu+1)}{2} (1-p)^2=1-p^\nu-r p^\nu(1-p)-\frac{\nu(\nu+1)}{2} (1-p)^2$$ is an increasing function of $p$ in $[0, 1]$ with value $0$ at $p=1$.
%Of course, when $\nu \in \N^+$ the inequality has an easy intuitive proof.

According to the results of Section \ref{FddConvergence}, the sum after the max in \eqref{JumpMaxBound} equals $\PP(X\ge2)$ where
\beq X\sim \begin{cases} NB\big(\frac{w_0}{k}+\mu, \frac{t_1+(b_0/k)}{t_2+(b_0/k)}\big) & \text{for Theorem \ref{PolyaPathRegime1}},\\
\text{Poisson}\big(\frac{2\delta}{b_0}\big) & \text{for Theorem \ref{PolyaPathRegime2}}, \\
NB\big(\frac{w_0}{k}+\mu, \frac{1-c^{-b_0-kt_1}}{1-c^{-b_0-kt_2}}\big) & \text{for Theorem \ref{QPolyaPathRegime1}},\\
\text{Poisson}\big(2\delta \frac{\log c}{c^{b_0}-1}\big) & \text{for Theorem \ref{QPolyaPathRegime2}}, \\
\end{cases}
\eeq
and $t_1:=(j-1)\delta, t_2:=(j+1)\delta$. The claim then follows from \eqref{PoissonTail} and \eqref{NBTail}.

\subsection{Conclusion}

It is clear from the form of the finite dimensional distributions that in all Theorems  \ref{PolyaPathRegime1}, \ref{PolyaPathRegime2},  \ref{QPolyaPathRegime1}, \ref{QPolyaPathRegime2} the limiting process $Z$ is a pure birth process that does not explode in finite time. Its rate at the point $(t, j)\in [0, \infty)\times \N$ is
$$\lambda_{t, j}=\lim_{h\to0^+}\frac{1}{h} \PP(Z(t+h)=j+1| Z(t)=j)$$
and is found as stated in the statement of each theorem.

\section{Deterministic and diffusion limits. Proof of Theorems \ref{PolyaPathLinear}, \ref{PolyaDiffusionLimit}, \ref{QPolyaPathODE}, \ref{QPolyaDiffusionLimit} } \label{DetermDiffLimits}

These theorems are proved with the use of Theorem 7.1 in Chapter 8 of \cite{Du96}, which is concerned with convergence of time-homogeneous Markov processes to diffusions. Since our basic Markov chain, $(A_n^{(m)})_{n\in\N}$, is not time-homogeneous, we do the standard trick of considering the chain  $\{(A_n^{(m)}, n)\}_{n\in\N}$ which is time-homogeneous.

\subsection{Proof of Theorems  \ref{PolyaPathLinear}, \ref{QPolyaPathODE} }

For each $m\in\N^+$, consider the discrete time-homogeneous Markov chain
$$Z_{n}^{(m)}=\Big(\frac{A_n^{(m)}}{m}, \frac{n}{m}\Big).$$
From any given state $(x_1, x_2)$ of $Z_n^{(m)}$, the chain moves to either of $(x_1+k/m, x_2+m^{-1}), (x_1, x_2+m^{-1})$ with corresponding probabilities $p(x_1, x_2, m), 1-p(x_1, x_2, m)$, where 
$$p(x_1, x_2, m):=\begin{cases} \frac{mx_1}{A_0^{(m)}+B_0^{(m)}+kmx_2} & \text{ in the case of Theorem \ref{PolyaPathLinear}},\\
\frac{1-q_m^{mx_1}}{1-q_m^{A_0^{(m)}+B_0^{(m)}+kmx_2}} & \text{ in the case of Theorem \ref{QPolyaPathODE}}.
\end{cases}
$$
This is true because when the chain is at the point $(x_1, x_2)$, then the time $n$ is $n=m x_2$ and $A_{n}^{(m)}+B_{n}^{(m)}=A_0^{(m)}+B_0^{(m)}+kn$. Define also
\begin{equation} p(x_1, x_2):=\lim_{m\to\infty} p(x_1, x_2, m)=\begin{cases} \frac{x_1}{a+b+kx_2} & \text{ in the case of Theorem \ref{PolyaPathLinear}},\\
\frac{1-c^{x_1}}{1-c^{a+b+kx_2}} & \text{ in the case of Theorem \ref{QPolyaPathODE}}.
\end{cases}
\end{equation}
We compute the mean and the covariance matrix for the one step change of $Z^{(m)}=(Z^{(m), 1}, Z^{(m), 2})$ conditioned on its current position.
\begin{align}
\EE\left[Z_{n+1}^{(m), 1}-Z_n^{(m), 1}|Z^{(m)}_n=\left(x_{1},x_{2}\right)\right]&=\frac{k}{m}p(x_1, x_2, m),\\ \EE\left[Z_{n+1}^{(m), 2}-Z_n^{(m), 2}|Z^{(m)}_n=\left(x_{1},x_{2}\right)\right]&=\frac{1}{m},\\
\EE\left[(Z_{n+1}^{(m), 1}-Z_n^{(m), 1})^2|Z^{(m)}_n=\left(x_{1},x_{2}\right)\right]&=\frac{k^2}{m^2}p(x_1, x_2, m),\\
\EE\left[(Z_{n+1}^{(m), 1}-Z_n^{(m), 1})(Z_{n+1}^{(m), 2}-Z_n^{(m), 2})|Z^{(m)}_n=\left(x_{1},x_{2}\right)\right]&=\frac{k}{m^{2}}p(x_1, x_2, m),\\
\EE\left[(Z_{n+1}^{(m), 2}-Z_n^{(m), 2})^2|Z^{(m)}_n=\left(x_{1},x_{2}\right)\right]&=\frac{1}{m^{2}}.
\end{align}
For each $m\in\Np$, we consider the process $\Lambda_t^{(m)}:=Z_{[mt]}^{(m)}, t\ge0$. According to Theorem 7.1 in Chapter 8 of \cite{Du96}, the sequence $(\Lambda^{(m)})_{m\ge1}$ converges weakly to the solution, $(S_{t})_{t\ge0}$, of the differential equation 
\begin{equation} \label{VectorODE}
\begin{aligned}
dS_{t}&=b(S_{t})dt,\\
S_0&=\Big(\begin{array}{c} a \\ 0\end{array}\Big),
\end{aligned}
\end{equation}
where
\begin{equation} 
S_{t}=\left(\begin{array}{c}S_{t}^{(1)} \\S_{t}^{(2)}\end{array}\right),  \,  b\Big(\begin{array}{c} x \\ y\end{array}\Big)=\Big(\begin{array}{c} kp(x, y) \\1\end{array}\Big).
\end{equation}
To apply the theorem, we need to check that the martingale problem MP$(b, \mathbb{O})$ has a unique solution. Here, $\mathbb{O}$ is the $2\times2$ zero matrix. See \cite{Du96}, \S 5.4, for details on the martingale problem. The problem indeed has unique solution because the differential equation \eqref{VectorODE} has a unique solution, and by well known results, this implies the claim for the martingale problem. 

It follows that the process $(A^{(m)}_{[mt]})_{t\ge0}$ converges, as $m\to\infty$, to the solution of the differential equation
\begin{align}X_{0}&=a,\\
dX_{t}&=k p(X_t, t)dt.
\end{align}
For both theorems, \ref{PolyaPathLinear} and \ref{QPolyaPathODE}, this ordinary differential equation is separable and its unique solution is the one stated.

\subsection{Proof of Theorems \ref{PolyaDiffusionLimit}, \ref{QPolyaDiffusionLimit}} 

\begin{proof}[\textbf{Proof of Theorem \ref{PolyaDiffusionLimit}}]

Call $\gl:=a/(a+b)$. 
For each $m\in\N^+$, consider the discrete time-homogeneous Markov chain
$$Z_{n}^{(m)}=\Big(\sqrt{m}\Big(\frac{A_n^{(m)}}{m}-a-\gl k \frac{n}{m}\Big), \frac{n}{m}\Big), n\in \N.$$

From any given state $(x_1, x_2)$ of $Z_n^{(m)}$, the chain moves to either of $(x_1-k  m^{-1/2}\gl, x_2+m^{-1}), (x_1+k m^{-1/2}(1-\gl), x_2+m^{-1})$ with corresponding probabilities 
\begin{align}
\Pi^{(m)}\left[\left(x_{1},x_{2}\right),\left(x_1-\frac{k}{\sqrt{m}}\gl, x_{2}+\frac{1}{m}\right)\right]&=\frac{B_n^{(m)}}{A_n^{(m)}+B_n^{(m)}},\\
\Pi^{(m)}\left[\left(x_{1},x_{2}\right),\left(x_1+\frac{k}{\sqrt{m}}(1-\gl), x_{2}+\frac{1}{m}\right)\right]&=\frac{A_n^{(m)}}{A_n^{(m)}+B_n^{(m)}},
\end{align}
with
\begin{align}
A_n^{(m)}&=ma+\gl k m x_2+x_1\sqrt{m},\\
B_n^{(m)}&=A_0^{(m)}+B_0^{(m)}+k m x_2-A_n^{(m)}.
\end{align}
We used the fact that when the chain is at the point $(x_1, x_2)$, then the time $n$ is $n=m x_2$.

We compute the mean and the covariance matrix for the one step change of $Z^{(m)}=(Z^{(m), 1}, Z^{(m), 2})$ conditioned on its current position.
\begin{align}
\EE\left[Z_{n+1}^{(m), 1}-Z_n^{(m), 1}|Z^{(m)}_n=\left(x_{1},x_{2}\right)\right]&=\frac{k}{\sqrt{m}} \frac{(1-\gl) A_n^{(m)}-\gl B_n^{(m)}}{A_n^{(m)}+B_n^{(m)}} \notag \\
 &\sim \frac{1}{m}\frac{k \{x_{1}-(\theta_1+\theta_2)\gl\}}{a+b+k x_{2}},\\
 \EE\left[Z_{n+1}^{(m), 2}-Z_n^{(m), 2}|Z^{(m)}_n=\left(x_{1},x_{2}\right)\right]&=\frac{1}{m},\\
\EE\left[(Z_{n+1}^{(m), 1}-Z_n^{(m), 1})^2|Z^{(m)}_n=\left(x_{1},x_{2}\right)\right]&=\frac{k^2}{m}\left(\gl^2 \frac{B_n^{(m)}}{A_n^{(m)}+B_n^{(m)}}+(1-\gl)^2 \frac{A_n^{(m)}}{A_n^{(m)}+B_n^{(m)}}\right) \notag \\
&\sim  \frac{1}{m} \frac{k^2 a b}{(a+b)^2},\\
\EE\left[(Z_{n+1}^{(m), 1}-Z_n^{(m), 1})(Z_{n+1}^{(m), 2}-Z_n^{(m), 2})|Z^{(m)}_n=\left(x_{1},x_{2}\right)\right]&\sim\frac{1}{m^{2}}\frac{kx_{1}}{a+b+kx_{2}},\\
\EE\left[(Z_{n+1}^{(m), 2}-Z_n^{(m), 2})^2|Z^{(m)}_n=\left(x_{1},x_{2}\right)\right]&=\frac{1}{m^{2}}.
\end{align}
Then, for each $m\in\N^+$, we consider the process $\Lambda_t^{(m)}:=Z_{[mt]}^{(m)}, t\ge0$. According to Theorem 7.1 in Chapter 8 of \cite{Du96}, the sequence $(\Lambda^{(m)})_{m\ge1}$ converges in distribution to the solution, $(S_{t})_{t\ge0}$, of the stochastic differential equation 
\begin{align}
dS_{t}&=b(S_{t})dt+\sigma(S_{t})dB_{t},\\
S_0&=\left(\begin{array}{c}\theta_1 \\0\end{array}\right),
\end{align}
where
$$\begin{array}{lll}
&S_{t}=\left(\begin{array}{c}S_{t}^{(1)} \\S_{t}^{(2)}\end{array}\right), & B_{t}=\left(\begin{array}{c}B_{t}^{(1)} \\B_{t}^{(2)}\end{array}\right),\\
&b\Big(\begin{array}{c} x \\ y\end{array}\Big)=\Big(\begin{array}{c} \frac{k \{x-(\theta_1+\theta_2)\gl\}}{a+b+k y} \\1\end{array}\Big), &  
\sigma\Big(\begin{array}{c} x \\ y\end{array}\Big)=\left(\begin{array}{ccc}
	k\frac{\sqrt{ab}}{a+b}&0\\
	0&0\\
\end{array}\right).
\end{array}$$
$B$ is a two dimensional standard Brownian motion.  Again, to apply that theorem, we need to check that the martingale problem MP$(b, \sigma)$ has a unique solution. This follows from the existence and uniqueness of strong solution for the above stochastic differential equation as the coefficients $b, \sigma$ are Lipschitz and grow at most linearly at infinity.

Thus, the process $(Z^{(m), 1}_{[mt]})_{t\ge0}$ converges in distribution, as $m\to\infty$, to the solution of
\begin{align}
Y_0&=\theta_1, \\
dY_t&=\frac{k\{Y_t-(\theta_1+\theta_2)\gl\}}{a+b+kt}\,dt+k\frac{\sqrt{ab}}{a+b}dB_t^{(1)}.
\end{align}
The same is true for $(C_t^{(m)})_{t\ge0}$ because $\sup_{t\ge0}|C_t^{(m)}-Z^{(m), 1}_{[mt]}|\le k/\sqrt{m}$. To solve the last SDE, we set $U_{t}:=\{Y_{t}-(\theta_1+\theta_2)\gl\}/(a+b+kt)$. 
Ito's lemma gives that
$$dU_{t}=k\frac{\sqrt{ab}}{(a+b)}\frac{1}{a+b+kt}dB_{t}^{(1)},$$
and since  $U_{0}=(b\theta_1-a\theta_2)/(a+b)^2$, we get
$$U_{t}=\frac{b\theta_1-a\theta_2}{(a+b)^2}+k\frac{\sqrt{ab}}{a+b} \int_{0}^{t}\frac{1}{a+b+ks}dB_{s}^{(1)}.$$
This gives \eqref{PolyaDiffLimit}. 
 \end{proof}

\begin{proof}[\textbf{Proof of Theorem \ref{QPolyaDiffusionLimit}}] The proof is analogous to that of Theorem \ref{PolyaDiffusionLimit}, only the algebra is a little more involved. For each $m\in\N^+$, consider the discrete time-homogeneous Markov chain
$$Z_{n}^{(m)}=\Big(\sqrt{m}\Big(\frac{A_n^{(m)}}{m}-X_{n/m}\Big), \frac{n}{m}\Big), n\in \N.$$
From any given state $(x_1, x_2)$ of $Z_n^{(m)}$, the chain moves to either of 
\begin{align}&(x_1, x_2)+(k  m^{-1/2}+\sqrt{m}(X_{x_2}-X_{x_2+m^{-1}}) , m^{-1}), \\
& (x_1, x_2)+(\sqrt{m}(X_{x_2}-X_{x_2+m^{-1}}) , m^{-1}) 
\end{align}
with corresponding probabilities $p(x_1, x_2, m), 1-p(x_1, x_2, m)$, where
\begin{equation}p(x_1, x_2, m)=
 \frac{[A_n^{(m)}]_{q_m}}{[A_0^{(m)}+B_0^{(m)}+kmx_2]_{q_m}} 
\end{equation}
and
\begin{align}
A_n^{(m)}&=mX_{x_2}+x_1\sqrt{m},\\
B_n^{(m)}&=A_0^{(m)}+B_0^{(m)}+k m x_2-A_n^{(m)}.
\end{align}
We used the fact that when the chain is at the point $(x_1, x_2)$, then the time $n$ is $n=m x_2$. For convenience, let $\Delta X_{x_2}=X_{x_2+m^{-1}}-X_{x_2}$.

We compute the mean and the covariance matrix for the one step change of $Z^{(m)}=(Z^{(m), 1}, Z^{(m), 2})$ conditioned on its current position. Of the following relations, the first four are immediate, the fifth follows from part (a) of the claim that follows and the fact that $Z_{n+1}^{(m), 2}-Z_n^{(m), 2}=1/m$.

\begin{align}
\EE\left[Z_{n+1}^{(m), 1}-Z_n^{(m), 1}|Z^{(m)}_n=\left(x_{1},x_{2}\right)\right]&=\frac{k}{\sqrt{m}} p(x_1, x_2, m)-\sqrt{m} \Delta X_{x_2} \label{Asymptotics1} \\
\EE\left[(Z_{n+1}^{(m), 1}-Z_n^{(m), 1})^2|Z^{(m)}_n=(x_1, x_2)\right]&=\left(\frac{k^2}{m}-2k \Delta X_{x_2}\right)p(x_1, x_2, m)+m (\Delta X_{x_2})^2   \label{Asymptotics2} \\
 \EE\left[Z_{n+1}^{(m), 2}-Z_n^{(m), 2}|Z^{(m)}_n=\left(x_{1},x_{2}\right)\right]&=\frac{1}{m},\\
\EE\left[(Z_{n+1}^{(m), 2}-Z_n^{(m), 2})^2|Z^{(m)}_n=\left(x_{1},x_{2}\right)\right]&=\frac{1}{m^{2}},\\
\EE\left[(Z_{n+1}^{(m), 1}-Z_n^{(m), 1})(Z_{n+1}^{(m), 2}-Z_n^{(m), 2})|Z^{(m)}_n=\left(x_{1},x_{2}\right)\right]&= O(m^{-2})
\end{align}
We examine now the asymptotics of the first two expectations.

\textsc{Claim}:
\begin{equation*} 
\begin{array}{lrl} (a) & \EE\left[Z_{n+1}^{(m), 1}-Z_n^{(m), 1}|Z^{(m)}_n=\left(x_{1},x_{2}\right)\right]\sim&
\dfrac{1}{m} \dfrac{k\log c}{c^{a+b+k x_2}-1}   \left(c^{X_{x_2}}x_1-\dfrac{(c^{X_{x_2}}-1)c^{a+b+k x_2}}{c^{a+b+k x_2}-1} (\theta_1+\theta_2)\right)\\ & & +O(\dfrac{1}{m^{3/2}})  \\ 
\phantom{adgagf}\\
 (b)  & \EE\left[\{Z_{n+1}^{(m), 1}-Z_n^{(m), 1}\}^2|Z^{(m)}_n=\left(x_{1},x_{2}\right)\right]
 \sim & \dfrac{1}{m} k^2 g(x_2)\{1-g(x_2)\}+O(\dfrac{1}{m^{3/2}}) 
\end{array} 
\end{equation*}
where $g(x_2):=\lim_{m\to\infty} p(x_1, x_2, m)=\frac{c^{X_{x_2}}-1}{c^{a+b+kx_2}-1}$.

\textsc{Proof of the claim}.
We examine the asymptotics of $p(x_1, x_2, m)$ and $\Delta X_{x_2}$. We have 
\begin{align} p(x_1, x_2, m)&=\frac{c^{X_{x_2}+\frac{1}{\sqrt{m}x_1}}-1}{c^{\frac{A_0^{(m)}+B_0^{(m)}}{m}+kx_2}-1}=\frac{c^{X_{x_2}+\frac{1}{\sqrt{m}x_1}}-1}{c^{a+b+kx_2+\frac{\theta_1+\theta_2}{\sqrt{m}}+O(\frac1m)}-1}\\
&=g(x_2)+\frac{\log c}{c^{a+b+k x_2}-1}   \left(c^{X_{x_2}}x_1-\frac{(c^{X_{x_2}}-1) c^{a+b+k x_2}}{c^{a+b+kx_2}-1}(\theta_1+\theta_2) \right)\frac{1}{\sqrt{m}}+O(\frac1m).
\end{align}
The second equality follows from a Taylor development. Also
\begin{equation}\Delta X_{x_2}=X'_{x_2}\frac{1}{m}+O(m^{-2})=k g(x_2)\frac{1}{m}+O(m^{-2}).
\end{equation}
For $X'$ we used the differential equation, \eqref{ODEqPolya}, that $X$ satisfies instead of the explicit expression for it. Substituting these expressions in \eqref{Asymptotics1}, \eqref{Asymptotics2}, we get the claim. 

Relation \eqref{qPolyaDetLimit} implies that $c^{X_{x_2}}=(c^{a+b}-1)/\{c^b-1+c^{-k x_2}(1-c^{-a})\}$, and this gives that the parenthesis following $\frac{1}{m}$ in equation (a) of the claim above equals 
\begin{equation}
\frac{(c^{a+b}-1) x_1-c^b(c^a-1)(\theta_1+\theta_2)}{c^b-1+c^{-kx_2}(1-c^{-a})}
\end{equation}
and also that
\begin{equation}g(x_2)\{1-g(x_2)\}=\frac{(c^a-1)(c^b-1) c^{a+kx_2}}{(c^{a+b+kx_2}-c^{a+kx_2}+c^a-1)^2}.
\end{equation}
It follows as before that the process $(Z^{(m)}_{[mt]})_{t\ge0}$ converges, as $m\to\infty$, to the solution of the stochastic differential equation 
\begin{align}
dS_{t}&=b(S_{t})dt+\sigma(S_{t})dB_{t},\\
S_0&=\left(\begin{array}{c}\theta_1 \\0\end{array}\right),
\end{align}
where
$$\begin{array}{lll}
&S_{t}=\left(\begin{array}{c}S_{t}^{(1)} \\S_{t}^{(2)}\end{array}\right), \, B_{t}=\left(\begin{array}{c}B_{t}^{(1)} \\B_{t}^{(2)}\end{array}\right), \phantom{} \\
&b\Big(\begin{array}{c} x \\ y\end{array}\Big)=\Big(\begin{array}{c} \frac{k\log c}{c^{a+b+k y}-1} \bigg\{\frac{(c^{a+b}-1) x-c^b(c^a-1)(\theta_1+\theta_2)}{c^b-1+c^{-ky}(1-c^{-a})}\bigg\}  \\1\end{array}\Big), &  \phantom{} \\
& \sigma\Big(\begin{array}{c} x \\ y\end{array}\Big)=\left(\begin{array}{ccc}
	k\sqrt{(c^a-1)(c^b-1)} \frac{c^{(a+ky)/2}}{c^{a+b+ky}-c^{a+ky}+c^a-1}&0\\
	0&0\\
\end{array}\right). &\phantom{}
\end{array}$$
$B$ is a two dimensional standard Brownian motion.  Again, the martingale problem MP$(b, \sigma)$ has a unique solution due to the form of the functions $b, \sigma$. And with analogous arguments as in Theorem \ref{PolyaDiffusionLimit}, we get that the process $(\hat C_t^{(m)})_{t\ge0}$ converges to the unique solution of the stochastic differential equation \eqref{QPolyaNoiseSDE}. To solve that, we remark that a solution of a stochastic differential equation of the form
\begin{equation} \label{LinearSDE} dY_t=(\alpha (t) Y_t+\beta(t))dt+\gamma(t) dW_t
\end{equation}
with $\alpha, \beta, \gamma:[0, \infty)\to\R$ continuous functions
is given by 
\begin{equation}
Y_t=e^{\int_0^t \alpha(s)\, ds}\left(Y_0+\int_0^t \beta(s)  e^{-\int_0^s \alpha (r)\, dr}\, ds+\int_0^t \gamma(s) e^{-\int_0^s \alpha (r)\, dr}\, dW_s\right).
\end{equation}
[To discover the formula, we apply It\'o's rule to $Y_t \exp\{-\int_0^t \alpha(s)\, ds\}$ and use \eqref{LinearSDE}.] Applying this formula for the values of $\alpha, \beta, \gamma$ dictated by \eqref{QPolyaNoiseSDE} we arrive at \eqref{QPolyaDiffLimit}.
\end{proof}

\section{Proofs for the $q$-P\'olya urn with many colors} \label{ProofMoreColors}

\begin{proof}[\textbf{Proof of Theorem \ref{PMFOfWeakColors}}] First, the equality of the expressions in \eqref{probmassforvectorAig}, \eqref{probmassforvector} follows from the definition of the $q$-multinomial coefficient. 

We will prove \eqref{probmassforvectorAig} by induction on $l$. When $l=2$, \eqref{probmassforvectorAig} holds because of \eqref{QPolyaWhiteDistr}. In that relation, we have $x_1=x,  x_2=n-x$. Assuming that \eqref{probmassforvectorAig} holds for $l\ge2$ we will prove the case $l+1$. The probability $\PP\left(X_{n,2}=x_{2},\ldots,X_{n,l+1}=x_{l+1}\right)$ equals 
\begin{align}
&\PP\left(X_{n,3}=x_3,\ldots,X_{n,l+1}=x_{l+1}\right)\PP(X_{n,2}=x_2 \mid X_{n,3}=x_3,\ldots,X_{n,l+1}=x_{l+1}) \label{PTimesConditional}\\
&=q^{\sum_{i=3}^{l+1}x_i\sum_{j=1}^{i-1}\left(a_{j}+kx_{j}\right)}\frac{ {-\frac{a_1+a_2}{k} \brack x_1+x_2}_{q^{-k}}  \prod_{i=3}^{l+1}{-\frac{a_{i}}{k} \brack x_{i}}_{q^{-k}}}{{-\frac{a_{1}+\ldots a_{l+1}}{k} \brack n}_{q^{-k}}}  \cdot q^{ x_2\left(a_1+k x_1\right)}\frac{ {-\frac{a_1}{k} \brack x_1}_{q^{-k}} {-\frac{a_2}{k} \brack x_2}_{q^{-k}}  }{{-\frac{a_1+a_2}{k} \brack x_1+x_2}_{q^{-k}}}
\\ &=q^{\sum_{i=2}^{l+1}x_{i}\sum_{j=1}^{i-1}\left(a_{j}+k x_{j}\right)}\frac{\prod_{i=1}^{l+1}{-\frac{a_{i}}{k} \brack x_{i}}_{q^{-k}}}{{-\frac{a_{1}+\ldots a_{l+1}}{k} \brack n}_{q^{-k}}}.
\end{align}
This finishes the induction provided that we can justify these two equalities. The second is obvious, so we turn to the first. The first probability in \eqref{PTimesConditional} is specified by the inductive hypothesis. That is, given the description of the experiment, in computing this probability it is as if we merge colors 1 and 2 into one color which is placed in the line before the remaining $l-1$  colors. This color has initially $a_1+a_2$ balls and we require that in the first $n$ drawings we choose it $x_1+x_2$ times.  The second probability in \eqref{PTimesConditional} is specified by the $l=2$ case of \eqref{probmassforvectorAig}, which we know. More specifically, since the number of drawings from colors $3, 4, \ldots, l+1$ is given, it is as if we have an urn with just two colors $1, 2$ that have initially $a_1$ and $a_2$ balls respectively. We do $x_1+x_2$ drawings with the usual rules for a $q$-P\'olya urn, placing in a line all balls of color 1 before all balls of color 2, and we want to pick $x_1$ times color 1 and $x_2$ times color 2.
\end{proof}

\begin{proof}[\textbf{Proof of Theorem \ref{LimitOfWeakColors}}]  The components of $(X_{n, 2}, X_{n, 3}, \ldots, X_{n, l})$ are increasing in $n$, and from Theorem \ref{LimitOfWeakColor} we have that each of them has finite limit (we treat all colors $2, \ldots, l$ as one color). Thus the convergence of the vector with probability one to a random vector with values is $\N^{l-1}$ follows. In particular, we also have convergence in distribution, and it remains to compute the distribution of the limit. Let $x_1:=n-(x_2+\cdots+x_l)$. Then the probability in \eqref{probmassforvectorAig} equals
\begin{align}
\PP\left(X_{n,2}=x_{2},\ldots,X_{n,l}=x_{l}\right)&=q^{-\sum_{1\le i<j \le l } a_j x_i }\frac{\prod_{i=1}^{l}{\frac{a_{i}}{k}+x_i-1 \brack x_{i}}_{q^{-k}}}{{\frac{\sum_{i=1}^{l}a_{i}}{k}+n-1 \brack n}_{q^{-k}}}\\&=q^{\sum_{1\le j< i \le l}x_{i} a_j}\frac{\prod_{i=1}^{l}{\frac{a_{i}}{k}+x_i-1 \brack x_{i}}_{q^{k}}}{{n+\frac{\sum_{i=1}^{l}a_{i}}{k}-1 \brack n}_{q^{k}}} \\
&=q^{\sum_{i=2}^{l}\left(x_{i}\sum_{j=1}^{i-1}a_{j}\right)}\left\{\prod_{i=2}^{l}{\frac{a_{i}}{k}+x_i-1 \brack x_{i} }_{q^{k}}\right\}\frac{{x_1+\frac{a_{1}}{k}-1 \brack x_{1}}_{q^{k}}}{{n+\frac{\sum_{i=1}^{l}a_{i}}{k}-1 \brack n}_{q^{k}}}. \label{FDDManyColors}
\end{align}
In the first equality, we used \eqref{NegToPos} while in the second we used \eqref{qMinusOneToq}. When we take $n\to\infty$ in \eqref{FDDManyColors}, the only terms involving $n$ are those of the last fraction, and \eqref{QBinomialAsymptotics} determines their limit.  Thus, the limit of  \eqref{FDDManyColors} is found to be the function $f(x_{2},\ldots,x_{l})$ in the statement of the theorem.
\end{proof}

% relation (2.2) in \cite{Char12} shows that $\sum_{x_{l}\in\N}\ldots\sum_{x_{2}\in\N} f(x_{2},\ldots,x_{l})=1$, so that $f$ %is a probabilty mass function of a random vector $\left(X_{2},\ldots,X_{l}\right)$ with values in $\N$.

\begin{proof}[\textbf{Proof of Theorem \ref{QPolyaPathODEWColors}}]
For each $m\in \Np$, we consider the discrete time-homogeneous Markov chain 
$$Z_n^{(m)}:=\left(\frac{n}{m}, \frac{A_{n,2}^{(m)}}{m}, \frac{A_{n,3}^{(m)}}{m}, \ldots, \frac{A_{n,l}^{(m)}}{m}\right), n\in \N.$$ 
From any given state $(t, x):=(t,x_2, x_3, \ldots, x_l)$ that $Z^{(m)}$ finds itself it moves to one of 
\begin{align*} &\left( t+\frac{1}{m}, x_2, \ldots, x_i+\frac{1}{m}, \ldots, x_l\right), i=2, \ldots, l,\\
&\left( t+\frac{1}{m}, x_2, \ldots, x_i, \ldots, x_l\right)
\end{align*}
with corresponding probabilities
\begin{align}
p_i(x_2, \ldots, x_l, t, m)&=q^{m s_{i-1}(t)}\frac{[m x_i]_q}{[m s_l(t)]_q},  i=2, \ldots, l,\\
p_1(x_2, \ldots, x_l, t, m)&=\frac{ [m x_1(t)]_q}{[m s_l(t)]_q},
\end{align}
where
$s_i(t)=x_1(t)+\sum_{1<j\le i} x_j$ for $i\in \{1, 2, \ldots, l\}$ and $x_1(t):=m^{-1}\sum_{j=1}^l A_{0, j}^{(m)}+kt-\sum_{2\le j \le l} x_i$. These follow from \eqref{problcolors} once we count the number of balls of each color present at the state $(t, x)$. To do this, we note that $Z_n^{(m)}=(t, x)$ implies that $n=mt$ drawings have taken place so far, the total number of balls is $A^{(m)}_{0, 1}+\cdots+A_{0, l}^{(m)}+k mt$, and the number of balls of color $i$, for $2\le i \le l$, is $m x_i$. Thus, the number of balls of color 1 is $A^{(m)}_{0, 1}+\cdots+A_{0, l}^{(m)}+k mt-m \sum_{2\le j\le l} x_i=mx_1(t)$. The required relations follow.

Let $x_1:=\lim_{m\to\infty} x_1(t)=\sigma_l+kt-\sum_{2\le j\le l} x_i$ and $s_i:=\lim_{m\to\infty} s_i(t)=\sum_{1\le j \le i}x_i$ for all $i\in  \{1, 2, \ldots, l\}$. Then, since $q=c^{1/m}$, for fixed $(t, x_2, \ldots, x_l)\in [0, \infty)^l$ with $(x_2, \ldots, x_l)\ne 0$, we have
\begin{equation}
\lim_{m\to\infty}p_i(x_2, \ldots, x_l, t, m)=c^{s_{i-1}} \frac{ [x_i]_c}{[s_l]_c}
\end{equation}
for all $i=2, \ldots, l$. We also note the following.
\begin{align}
Z_{{n+1}, 1}^{(m)}-Z_{n, 1}^{(m)}&=\frac{1}{m},\\
\EE\left[Z_{{n+1}, i}^{(m)}-Z_{n, i}^{(m)}|Z_n^{(m)}=\left(t, x_{2},\ldots,x_{l}\right)\right]&=\frac{k}{m} p_i(x_2, \ldots, x_l, t, m),\\
\EE\left[(Z_{{n+1}, i}^{(m)}-Z_{n, i}^{(m)})^2|Z_n^{(m)}=\left(t, x_{2},\ldots,x_{l}\right)\right]&=\frac{k^2}{m^2} p_i(x_2, \ldots, x_l, t, m),\\
\EE\left[(Z_{{n+1}, i}^{(m)}-Z_{n, i}^{(m)})(Z_{{n+1}, j}^{(m)}-Z_{n, j}^{(m)})|Z_n^{(m)}=\left(t, x_{2},\ldots,x_{l}\right)\right]&=0
\end{align}
for $i,j=2,3,\ldots,l$ with $i\neq j$.

Therefore, with similar arguments as in the proof of Theorem \ref{PolyaPathLinear}, as $m\rightarrow+\infty , (Z_{[mt]}^{(m)})_{t\ge0}$ converges in distribution to $Y$, the solution of the ordinary differential equation
\begin{equation} \label{ODEColors} \begin{aligned}
dY_{t}&=b(Y_{t})dt,\\
Y_0&=(0, a_2, \ldots, a_l),
\end{aligned}
\end{equation}
where $b(t, x_2, \ldots, x_l)=\left(1, b^{(2)}(t, x),b^{(3)}(t, x),\ldots,b^{(l)}(t, x)\right)$ with
$$b^{(i)}(t, x)=k c^{s_{i-1}} \frac{ [x_i]_c}{[s_l]_c}$$
for $i=2,3,\ldots,l.$ Note that $s_l=\sigma_l+kt$ does not depend on $x$.
 
 Since $A^{(m)}_{[mt], 1}+A^{(m)}_{[mt], 2}+\cdots+A^{(m)}_{[mt], l}=kmt+A^{(m)}_{0, 1}+A^{(m)}_{0, 2}+\cdots+A^{(m)}_{0, l}$, we get that the process $(A^{(m)}_{[mt], 1}/m, A^{(m)}_{[mt], 2}/m+\cdots+A^{(m)}_{[mt], l}/m)_{t\ge0}$ converges in distribution to a process $(X_{t, 1}, X_{t, 2}, \ldots, X_{t, l})_{t\ge0}$ so that $X_{t, 1}+\cdots+X_{t, l}=a_1+a_2+\cdots+a_l+kt$, while the $X_{t, i}, i=2, \ldots, l$, satisfy the system
\begin{align}
X_{t, i}'&=kc^{\sigma_l+k t-\sum_{j=i}^l X_{t,i}}\frac{1-c^{X_{t,i}}}{1-c^{\sigma_l+kt}} \qquad \text{ for all } t>0,\\
X_{0, i}&=a_i,
\end{align}
with $i=2 ,3,\ldots,l$. Letting $Z_{r, i}=c^{X_{\frac{1}{k\log c}\log r, i}}$ for all $r\in(0, 1]$ and $i\in\{1, 2, \ldots, l\}$, we have for the $Z_{r, i}, i\in\{2, 3, \ldots, l\}$ the system
\begin{align}
\frac{Z'_{r, i}}{1-Z_{r, i}}&=\frac{\sigma_l}{1-\sigma_l r}\frac{1}{\prod_{i<j\le l} Z_{r, j}},\\
Z_{1, i}&=c^{a_i}.
\end{align}
In the case $i=l$, the empty product equals 1. It is now easy to prove by induction (starting from $i=l$ and going down to $i=2$) that 
\begin{equation} \label{ZSystemSol}
Z_{r, i}=\frac{c^{\sigma_l-\sigma_{i-1}}(1-c^{\sigma_l} r)-c^{\sigma_l}(1-r)}{c^{\sigma_l-\sigma_i} (1-c^{\sigma_l}r)-c^{\sigma_l}(1-r)}
\end{equation}
for all $r\in(0, 1]$. Since $Z_{r, 1}Z_{r, 2}\cdots Z_{r, l}=c^{\sigma_l} r$, we can check that \eqref{ZSystemSol} holds for $i=1$ too.
The fraction in  \eqref{ZSystemSol} equals
\begin{equation}c^{a_i} \frac{(1-c^{\sigma_l} r)-c^{\sigma_{i-1}}(1-r)}{(1-c^{\sigma_l}r)-c^{\sigma_i}(1-r)}.
\end{equation}
Recalling that $X_{t, i}=(\log c)^{-1} \log Z_{c^{kt}}$, we get \eqref{SystemSolution} for all $i\in\{1, 2, \ldots, l\}$ . 
\end{proof}

\begin{proof}[\textbf{Proof of Theorem \ref{QPolyaPathODEWColors2}}] This is proved in the same way as Theorem  \ref{QPolyaPathODEWColors}. We keep the same notation as there. The only difference now is that $\lim_{m\to\infty} p_i(t, x_2, \ldots, x_l, m)=x_i/s_l$. As a consequence, the system of ordinary differential equations for the limit process $Y_t:=(t, X_{t, 2}, \ldots, X_{t, l})$ is \eqref{ODEColors} but with
$$b^{(i)}(t, x)=\frac{kx_i}{s_l}.$$
Recall that $s_l=\sigma_l+kt$. 
Thus, for $i=2, 3, \ldots, l$, the process $X_{t, i}$ satisfies $X_{t, i}'=k X_{t, i}/(\sigma_l+kt), X_{0, i}=a_i$, which give immediately the last $l-1$ coordinates of \eqref{SystemSolutionReg2}. The formula for the first coordinate follows from $X_{t, 1}+X_{t, 2}+\cdots+X_{t, l}=kt+\sigma_l$.
\end{proof}


\begin{thebibliography}{99}

\bibitem{BaPa85} Bagchi, Arunabha, and Asim K. Pal. ``Asymptotic normality in the generalized P\'olya-Eggenberger urn model, with an application to computer data structures.'' SIAM Journal on Algebraic Discrete Methods 6, no. 3 (1985): 394-405.

\bibitem{BaHu02} Bai, Zhi-Dong, Feifang Hu, and Li-Xin Zhang. ``Gaussian approximation theorems for urn models and their applications." The Annals of Applied Probability 12, no. 4 (2002): 1149-1173.


\bibitem{Char12} Charalambides, Ch A. ``A q-P\'olya urn model and the $q$-P\'olya and inverse q-P\'olya distributions.'' Journal of Statistical Planning and Inference 142, no. 1 (2012): 276-288.

\bibitem{Char16} Charalambides, Charalambos A. Discrete $q$-distributions. John Wiley \& Sons, 2016.

\bibitem{Du96} Durrett, Richard. Stochastic calculus: a practical introduction. CRC press, 1996.

\bibitem {EtKu} Ethier, Stewart N., and Thomas G. Kurtz. Markov processes: characterization and convergence. Vol. 282. John Wiley \& Sons, 2009.
\bibitem {Gas} G. Gasper and M. Rahman. Basic Hypergeometric Series. Cambridge University Press, Cambridge, 1990

\bibitem {Go93} Gouet, Raul. ``Martingale functional central limit theorems for a generalized P\'olya urn." The Annals of Probability (1993): 1624-1639.

\bibitem{He77} Heyde, C. C. ``On central limit and iterated logarithm supplements to the martingale convergence theorem." Journal of Applied Probability 14, no. 4 (1977): 758-775.

\bibitem{Ja} Janson, Svante. ``Functional limit theorems for multitype branching processes and generalized P\'olya urns." Stochastic Processes and their Applications 110, no. 2 (2004): 177-245.

\bibitem{JoKo} Johnson, Norman Lloyd, and Samuel Kotz. ``Urn models and their application; an approach to modern discrete probability theory." (1977).

\bibitem{KacCheung} Kac, Victor, and Pokman Cheung. Quantum calculus. Springer Science \& Business Media, 2001.

\bibitem{Kup} Kupershmidt, Boris A.  ``$q$-probability: I. Basic discrete distributions.'' Journal of Nonlinear Mathematical Physics 7, no. 1 (2000): 73-93.

\bibitem{Mah} Mahmoud, Hosam. P\'olya urn models. Chapman and Hall/CRC, 2008.


\end{thebibliography}
\end{document}